\documentstyle[12pt,amsmath,amsfonts,amssymb]{article}
\newtheorem{theorem}{Theorem}[section]
\newtheorem{definition}[theorem]{Definition}
\newtheorem{corollary}[theorem]{Corollary}
\newtheorem{proposition}[theorem]{Proposition}
\newtheorem{remark}[theorem]{Remark}
\newtheorem{lemma}[theorem]{Lemma}

\newcommand {\Uc}      {{\mathcal U}}
\newcommand {\Vc}      {{\mathcal V}}

\newcommand {\Hc}      {{\mathcal H}}
\newcommand {\Ac}      {{\mathcal A}}
\newcommand {\Tc}      {{\mathcal T}}

\newcommand {\Pc}      {{\mathcal P}}
\newcommand {\Ec}      {{\mathcal E}}
\newcommand {\Qc}      {{\mathcal Q}}
\newcommand {\WKP}     {W^k_p(\RN)}
\newcommand {\FS}      {F^s_{pq}(\RN)}
\newcommand {\BS}      {B^s_{pq}(\RN)}

\newcommand {\PK}      {\Pc_k}

\newcommand {\tK}      {\widetilde{K}}
\newcommand {\tQ}      {\widetilde{Q}}

\newcommand {\tf}      {\tilde{f}}

\newcommand {\R}       {{\bf R}}
\newcommand {\N}       {{\bf N}}
\newcommand {\RN}      {\R^n}
\newcommand {\fs}      {f\,^\sharp_{\vv,S}}
\newcommand {\fsd}     {f\,^\sharp_{\vv,\Delta,S}}
\newcommand {\tfsd}    {(\tf)\,^\sharp_{\vv,\Delta,\RN}}
\newcommand {\tfs}     {(\tf)^\sharp_{\vv}}
\newcommand {\1}       {\dist(x,S)}
\newcommand {\ve}      {\varepsilon}
\newcommand {\vv}      {\vec{v}}
\newcommand {\emp}     {\emptyset}
\newcommand {\cw}      {\curlywedge}

\newcommand {\QB}      {\overline{Q}_K}
\newcommand {\OK}      {\overline{K}}
\newcommand {\ip}[1]   {\langle{#1}\rangle}
\newcommand {\card}    {\operatorname{card}}
\newcommand {\supp}    {\operatorname{supp}}
\newcommand {\Ex}      {\operatorname{Ext}}

\newcommand {\Eks}     {\Ex_{k,S}}
\newcommand {\Ekus}    {\Ex_{k,u,\,S}}
\newcommand {\diam}    {\operatorname{diam}}
\newcommand {\dist}    {\operatorname{dist}}
\newcommand {\cl}      {\operatorname{cl}}
\newcommand {\PR}      {\operatorname{Pr}}
\newcommand {\bx}      {\hspace{10mm}$\Box$}
\newcommand {\BX}      {\hspace{10mm}\Box}
\newcommand {\nn}      {\nonumber}
\newcommand {\rf}[1]   {(\ref{#1})}
\newcommand {\SECT}[2] {\section*{\centerline{\normalsize
{\bf #1}}} \setcounter{section}{#2}
\setcounter{theorem}{0}\setcounter{equation}{0}}
\newcommand {\SECTLONG}[3]
{\section*{\centerline{\normalsize {\bf #1}}
\centerline{\normalsize {\bf #2}}} \setcounter{section}{#3}
\setcounter{theorem}{0}\setcounter{equation}{0}}
\newcommand{\be}            {\begin{eqnarray}}
\newcommand{\bel}[1]        {\begin{eqnarray}\label{#1}}
\newcommand{\ee}            {\end{eqnarray}}
\begin{document}
\medskip
\centerline{\large{\bf Local approximations and intrinsic
characterizations }} \vspace*{4mm}
\centerline{\large{\bf of spaces of smooth functions on
regular subsets of $\RN$}}
\vspace*{10mm} \centerline{By  {\it Pavel Shvartsman}}
\vspace*{3mm}
\centerline {\it Department of Mathematics, Technion -
Israel Institute of Technology,}
\centerline{\it 32000 Haifa, Israel}\vspace*{3 mm}
\centerline{\it e-mail: pshv@tx.technion.ac.il}
\vspace*{12 mm}
\renewcommand{\thefootnote}{ }
\footnotetext[1]{{\it\hspace{-6mm}Math Subject
Classification} 46E35}
\begin{abstract}
{\small
\par We give an intrinsic characterization of the
restrictions of Sobolev $\WKP$, Triebel-Lizorkin $\FS$ and
Besov $\BS$ spaces to regular subsets of $\RN$ via sharp
maximal functions and local approximations.}
\end{abstract}
\vspace*{15mm}
\renewcommand{\thefootnote}{\arabic{footnote}}
\setcounter{footnote}{0}
\SECT{1. Main definitions and results.}{1}
\indent
\par The purpose of this paper is to study the problem of
extendability of functions defined on measurable subsets of
$\RN$ to functions defined on the whole space and
satisfying certain smoothness conditions.
\par We will consider three kinds of spaces of smooth functions
defined on $\RN$. First we deal with classical Sobolev
spaces, see e.g. Maz'ja \cite{M}. We recall that, given an
open set $\Omega\subset\RN,$ $k\in\N$ and $p\in[1,\infty]$,
the Sobolev space $W^k_p(\Omega)$ consists of all functions
$f\in L_{1,\,loc}(\Omega)$ whose distributional partial
derivatives on $\Omega$ of all orders up to $k$ belong to
$L_p(\Omega)$. $W^k_p(\Omega)$ is normed by
$ \|f\|_{W^k_p(\Omega)}:=
 \sum\{\|D^\alpha f\|_{L_p(\Omega)}:|\alpha|\le k\}. $
\par There is an extensive literature devoted to
describing the restrictions of Sobolev functions to
different classes of subsets of $\RN$. (We refer the reader
to \cite{M,MP,BIN,BK,Bur,FJ,Ch,R1} and references therein
for numerous results and techniques in this direction.) Let
us recall some of these results. Calder\'{o}n \cite{C2}
showed that, if $\Omega$ is a Lipschitz domain in $\RN$ and
$1<p<\infty$, then $\WKP|_\Omega=W^k_p(\Omega).$ Stein
\cite{St} extended this result for $p=1,\infty$ and Jones
\cite{Jn} showed that the same isomorphism also holds for
every $(\ve,\delta)-$domain and $1\le p\le \infty$.
\par Here as usual, for any Banach space
$(\Ac, \|\cdot\|_{\Ac})$ of measurable functions defined on
$\RN$ and a measurable set $S\subset \RN$ of positive
Lebesgue measure, we let $\Ac|_S$ denote the restriction of
$\Ac$ to $S$, i.e., the Banach space
$$
\Ac|_S:=\{f:S\to\R~: {\rm there~is}\ \  F\in \Ac \ \ {\rm
such~that}\ \  F|_S=f\ \  \rm{~~a.~e.}\}.
$$
We call $\Ac|_S$ the {\it restriction space} or the {\it
trace space} of $\Ac$ to $S$. $\Ac|_S$ is equipped with the
standard quotient space norm
$$
\|f\|_{\Ac|_S}:=\inf\{\|F\|_{\Ac}: F\in
\Ac,F|_S=f\rm{~~a.~e.}\}.
$$
\par Our aim is to extend these results to the case of
so-called {\it regular subsets of $\RN$}. We define these
sets as follows.
\begin{definition}\label{Reg}
A measurable set $S\subset\RN$ is said to be regular if
there are constants $\theta_S\ge 1$ and $\delta_S>0$ such
that, for every cube $Q$ with center in $S$ and with
diameter
 $\diam Q\le\delta_S$,
$$
|Q|\le \theta_S |Q\cap S|.
$$
\end{definition}
\par Here $|A|$ stands for the Lebesgue measure of a
measurable set $A\subset\RN$. We will also assume that all
cubes in this paper are closed and have sides which are
parallel to the coordinate axes. We let $Q(x,r)$ denote the
cube in $\RN$ centered at $x$ with side length $2r$.
\par Regular subsets
of $\RN$ are often called Ahlfors $n$-regular or $n$-sets
\cite{JW}. Cantor-like sets and Sierpi\'{n}ski type gaskets
(or carpets) of positive Lebes\-gue measure are examples of
non-trivial regular subsets of $\RN$ .
\par We observe that the interior of a regular set can be
empty (as, for instance, for a Cantor-like set or for a
Sierpi\'{n}ski type gasket). Thus, to give a constructive
characterization of the restrictions of Sobolev functions
to regular sets, we need an equ\-i\-va\-lent definition of
the Sobolev spaces which does not use the notion of
derivatives.
\par There are several known ways of defining Sobolev
spaces which do not use derivatives. In this paper our
point of departure will be a characterization of Sobolev
spaces due to Calder\'on. In \cite{C} (see also \cite{CS})
Calder\'on characterizes the Sobolev spaces $\WKP$ in terms
of $L_p$-properties of sharp maximal functions.
\par
Before we recall Calder\'on's result we need to fix some
notation. Let $\PK=\Pc_k(\RN),$ $k\ge 0,$ denote the family
of all polynomials on $\RN$ of degree less than or equal to
$k$. We also put $\Pc_{-1}:=\{0\}$. Given $f\in
L_{u,\,loc}(\RN)$, $0<u\le\infty,$ and a cube $Q$, we let
$\Ec_k(f;Q)_{L_u}$ denote the {\it normalized local best
approximation} of $f$ on $Q$ in $L_u$-norm by polynomials
of degree at most $k-1$, see Brudnyi \cite{Br1}. More
explicitly, we define
\bel{ESP}
\Ec_k(f;Q)_{L_u}:=|Q|^{-\frac{1}{u}}\inf_{P\in\Pc_{k-1}}
\|f-P\|_{L_u(Q)}=\inf_{P\in\Pc_{k-1}}\left(\frac{1}{|Q|}
\int_Q|f-P|^u\,dx\right)^{\frac{1}{u}}. \ee

\par In the literature $\Ec_k(f;Q)_{L_u}$ is also sometimes
called the {\it local oscillation} of $f$, see e.g. Triebel
\cite{T}. This quantity is the main object of the theory of
{\it local polynomial approximations} which provides a
unified framework for the description of a large family of
spaces of smooth functions. We refer the reader to Brudnyi
\cite{Br3}-\cite{Br2} for the main ideas and results in
local approximation theory; see also Section 5 for
additional information and remarks related to this theory.
\par Given $\alpha>0$ and a locally integrable function $f$
on $\RN$, we define its {\it fractional sharp maximal
function} $f^{\sharp}_\alpha$ by letting
\bel{SM}
f^{\sharp}_\alpha(x):=\sup_{r>0}
r^{-\alpha}\,\Ec_{k}(f;Q(x,r))_{L_1}.
\ee
Here $k:=-[-\alpha]$ is the greatest integer strictly less
than $\alpha+1$.
\par In \cite{C} Calder\'{o}n proved
that, for $1<p\le \infty$, a function $f$ is in $\WKP$, if
and only if $f$ and $f^{\sharp}_k$ are both in $L_p(\RN)$.
Moreover, up to constants depending only on $n,k$ and $p$
the following equivalence,
\bel{AN}
\|f\|_{\WKP}\approx
\|f\|_{L_p(\RN)}+\|f_{k}^{\sharp}\|_{L_p(\RN)},
\ee
holds.
\par This characterization motivates the following
definition. Given $u>0$, a function $f\in L_{u,\,loc}(S),$
and a cube $Q$ whose center is in $S$, we let
$\Ec_k(f;Q)_{L_u(S)}$ denote the norma\-lized best
approximation of $f$ on $Q$ in $L_u(S)$-norm:
\bel{ES}
\Ec_k(f;Q)_{L_u(S)}:=|Q|^{-\frac{1}{u}}
\inf_{P\in\PK}\|f-P\|_{L_u(Q\cap S)}
=\inf_{P\in\Pc_{k-1}}\left(\frac{1}{|Q|} \int_{Q\cap S}
|f-P|^u\,dx\right)^{\frac{1}{u}}.
\ee
By $f^{\sharp}_{\alpha,S},$ we denote the fractional sharp
maximal function of $f$ on $S$,
$$
f^{\sharp}_{\alpha,S}(x):=\sup_{r>0}
r^{-\alpha}\,\Ec_{k}(f;Q(x,r))_{L_1(S)},\ \ \ \ \ \ x\in S.
$$
Here $k(=-[-\alpha])$ is the same as in \rf{SM}. Thus,
$f^{\sharp}_{\alpha}=f^{\sharp}_{\alpha,\RN}$~.
\par The first main result of the paper is the following
\begin{theorem}\label{EXT1}
Let $S$ be a regular subset of $\RN$. Then a function $f\in
L_p(S),~1<p\le \infty,$ can be extended to a function $F\in
\WKP$ if and only if
$$
f_{k,S}^{\sharp}:=\sup_{r>0}
r^{-k}\,\Ec_{k}(f;Q(\cdot,r))_{L_1(S)}\in L_p(S).
$$
In addition,
\bel{Sob}
\|f\|_{\WKP|_S}\approx
\|f\|_{L_p(S)}+\|f_{k,S}^{\sharp}\|_{L_p(S)}
\ee
with constants of equivalence depending only on
$n,k,p,\theta_S$ and $\delta_S$.
\end{theorem}
\par
For $k=1$ this theorem follows from a more general result
proved in \cite{S3} for the case of a metric space equipped
with a doubling measure.
\par We now turn to the second kind of
spaces of smooth functions to be considered in this paper,
namely the Triebel-Lizorkin spaces $\FS$. The reader can
find a detailed treatment of the theory of these spaces in
the monographs \cite{T,T3,ET,RS}. The scale $\FS$ includes,
in particular, the Bessel potential spaces
$H^s_p(\RN)=F^s_{p,2},1<p<\infty,$ (\cite{T}, p. 11). These
spaces which are also referred to in the literature as {\it
fractional Sobolev spaces} are generalizations of the
Sobolev spaces in the following sense: $H^k_p(\RN)=\WKP$
whenever $k\in\N$ and $1<p<\infty$.
\par Among the various equivalent definitions
of Triebel-Lizorkin spaces, the most useful one for us is
expressed in terms of local polynomial approximations:
\par Given $0<s<k$,
$1<p<\infty$, $1\le q\le\infty,$ a function $f\in
L_{1,\,loc}(\RN)$ and $x\in\RN$, we put
$$
g(x):=\left(\int_0^1
\left(\frac{\Ec_k(f;Q(x,t))_{L_1}}{t^s}\right)^q\,
\frac{dt}{t}\right)^{\frac{1}{q}}
$$
for $q<\infty$ and
$ g(x):=\sup\{t^{-s}\,\Ec_k(f;Q(x,t))_{L_1}:~t>0\}$
for $q=\infty$. Then $f\in \FS$ if and only if $f$ and $g$
are both in $L^p(\RN)$. Moreover,
\bel{Crit} \|f\|_{F^s_{pq}(\RN)}\approx
\|f\|_{L_p(\RN)}+\|g\|_{L_p(\RN)}
\ee
with constants depending only on $n,s,p,q$ and $k$. This
description is due to Dorronsoro \cite{Dor1,Dor2}, Seeger
\cite{See} and Triebel \cite{T2}; see also \cite{T}, p. 51,
and references therein for a detailed history of the
problem.
\par The second main result of the paper, Theorem
\ref{EXT2}, states that for a {\it regular} subset
$S\subset\RN$, the trace space $\FS|_S$ can be
characterized in a way which is analogous to the preceding
definition , i.e., in terms of local approximations of
functions taken on the set $S$ instead of on $\RN$.
\begin{theorem}\label{EXT2}
Let $S$ be a regular subset of $\RN$, $0<s<k$,
$1<p<\infty$, and $1\le q\le\infty$. Then a function $f\in
L_p(S)$ can be extended to a function $F\in\FS$ if and only
if
$$
\left(\int_0^1
\left(\frac{\Ec_k(f;Q(\cdot,t))_{L_1(S)}}{t^s}\right)^q\,
\frac{dt}{t}\right)^{\frac{1}{q}}\in L_p(S)
$$
(with usual modification for $q=\infty$). In addition,
\bel{FSnorm} \|f\|_{\FS|_S}\approx \|f\|_{L^p(S)}+
\left\|\left(\int_0^1
\left(\frac{\Ec_k(f;Q(\cdot,t))_{L_1(S)}}{t^s}\right)^q\,
\frac{dt}{t}\right)^{\frac{1}{q}}\right\|_{L_p(S)} \ee
with constants of equivalence depending only on
$n,\theta_S,\delta_S,s,p,q$ and $k$.
\end{theorem}
\par Observe that for Lipschitz domains in
$\RN$ an intrinsic characterization of traces of $F$-spaces
was given by Kalyabin \cite{Kal}. For
$(\ve,\delta)$-domains such a characterization is due to
Seeger \cite{See}; see also Triebel \cite{T2}.
\par Let us note two particular cases of Theorem
\ref{EXT2}.
\begin{remark} Recall that $\WKP=F^k_{p,\,2},$
$1<p<\infty$. This and Theorem \ref{EXT2} imply one more
intrinsic characterization of the restrictions of Sobolev
functions to regular subsets, cf. \rf{Sob}: for every
$1<p<\infty$
$$
\WKP|_S=\{f\in L_p(S):~ \left(\int_0^1
\left(\frac{\Ec_k(f;Q(\cdot,t))_{L_1(S)}}{t^s}\right)^2\,
\frac{dt}{t}\right)^{\frac{1}{2}}\in L_p(S)\}.
$$
\end{remark}
\begin{remark} The space $F^s_{p,\infty}(\RN),$
$s>0,1\le p\le\infty$, coincides with the space
$C^s_p(\RN)$ introduced by DeVore and Sharpley \cite{DS};
for non-integer $s$ this space was independently
consi\-de\-red by Christ \cite{Ch}. (See also Triebel
\cite{T}, p. 48--50 for additional remarks and comments.)
We recall that $C^s_p(\RN)$ consists of all functions $f$
defined on $\RN$ such that $f, f^\flat_s\in L^p(\RN)$. Here
$f^\flat_s$ is defined by formula \rf{SM} with $k=[s]+1$.
\par Thus by Theorem \ref{EXT2} for every regular set $S$,
every $0<s<k$ and $1< p<\infty$,
$$
C^s_p(\RN)|_S=\{f\in L^p(S):~ \sup_{0<t\le 1}
t^{-s}\Ec_k(f;Q(\cdot,t))_{L_1(S)}\in L_p(S)\}.
$$
Observe that Devore and Sharpley \cite{DS} obtained an
intrinsic characterization of the trace space
$C^s_p(\RN)|_\Omega$ where $\Omega$ is a Lipschitz domain
in $\RN$; for $(\varepsilon,\delta)$-domain it was done by
Christ \cite{Ch}.
\end{remark}
\par Having dealt with Sobolev and Triebel-Lizorkin spaces
we now turn finally to consider Besov spaces $\BS$. For a
general theory of these spaces we refer the reader to the
monographs \cite{BIN,T,RS} and references therein. See also
Section 5 for definitions and a description of the Besov
spaces via local approximations. This description provides
the following equivalent norm on the Besov spaces: for all
$1\le u\le p\le\infty$, $0<q\le \infty,$ and $0<s<k$
\bel{ChB} \|f\|_{\BS}\approx
\|f\|_{L_p(\RN)}+\left(\int_0^1
\left(\frac{\|\Ec_k(f;Q(\cdot,t))_{L_u}\|_{L_p(\RN)}}
{t^s}\right)^q\, \frac{dt}{t}\right)^{\frac{1}{q}} \ee
(with usual modification if $q=\infty$). Here constants of
the equivalence depend only on $n,s,k,p$ and $q$. This
characterization (in an equivalent form, via so-called
$(k,p)$-modulus of continuity in $L_u$ , see \rf{KPU}) was
obtained by Brudnyi \cite{Br1}; we also refer to Triebel
\cite{T}, p. 51, and references therein.
\par Our next result, Theorem \ref{EXT3}, states that,
similar to Sobolev and $F$-spaces, a natural generalization
of description \rf{ChB} to regular sets provides an
intrinsic characterization of the restrictions of Besov
functions.
\begin{theorem}\label{EXT3}
Let $S$ be a regular subset of $\RN$, $0<s<k$, $1\le u\le
p\le\infty$ and $0< q\le\infty$. Then a function
 $f\in L_p(S)$ can be extended to a function $F\in\BS$ if
and only if
$$
\int_0^1
\left(\frac{\|\Ec_k(f;Q(\cdot,t))_{L_u(S)}\|_{L_p(S)}}
{t^s}\right)^q\, \frac{dt}{t}<\infty
$$
($\sup\limits_{0<t\le
1}t^{-s}\|\Ec_k(f;Q(\cdot,t))_{L_u(S)}\|_{L_p(S)}<\infty$
for $q=\infty$). In addition,
$$
\|f\|_{\BS|_S}\approx \|f\|_{L^p(S)}+ \left(\int_0^1
\left(\frac{\|\Ec_k(f;Q(\cdot,t))_{L_u(S)}\|_{L_p(S)}}
{t^s}\right)^q\, \frac{dt}{t}\right)^{\frac{1}{q}}
$$
(modification if $q=\infty$). Here constants of equivalence
depend only on $n,\theta_S,\delta_S,s,p,q$ and $k$.
\end{theorem}
\par For intrinsic description of the Besov spaces on
Lipschitz domains we refer the reader to Nikol'ski \cite{N}
and Besov \cite{B1,B2}; see also Rychkov \cite{R2}. The
case of $(\ve,\delta)$-domains was independently treated by
Shvartsman \cite{S2}, Seeger \cite{See} and Devore and
Sharpley \cite{DS2}.
\par Proofs of Theorems \ref{EXT1}, \ref{EXT2} and
\ref{EXT3} are based on a modification of the Whitney
extension method suggested in author's work \cite{S1}, see
also \cite{S2}.
\par A crucial step of this approach is presented in
Section 2. Without loss of generality we may assume that
$S$ is closed so that $\RN\setminus S$ is open. By
$W_S=\{Q_k\}$ we denote a Whitney decomposition of
$\RN\setminus S$, see e.g. \cite{St}.
\par We assign every cube $Q=Q(x_Q,r_Q)\in W_S$
a measurable subset $H_Q\subset S$ such that $H_Q\subset
Q(x_Q,10r_Q)\cap S$, $|Q|\le \gamma_1|H_Q|$ whenever $\diam
Q\le \delta_S$, and a family of sets
$ \Hc_S:=\{H_Q:Q\in W_S\} $
has {\it a finite covering multiplicity }, i.e., every
point $x\in S$ belongs to at most $\gamma_2$ sets of the
family $\Hc_S$. Here $\gamma_1,\gamma_2$ are positive
constants depending only on $n$ and $\theta_S$. We call
every set $H_Q\in\Hc_S$ a {\it ``reflected quasi-cube"}
associated to the Whitney cube $Q$. The existence of the
family $\Hc_S$ of reflected quasi-cubes is proved in
Theorem \ref{HB}.
\par The second step of the extension method
is presented in Section 3. We associate every function
$f\in L_{1,\,loc}(S),$ and every Whitney cube $Q\in W_S$ of
$\diam Q\le\delta_S,$ a {\it linear} mapping
$P_Q:L_1(H_Q)\to \Pc_{k-1}$ which provides almost the best
polynomial approximation of $f$ on $H_Q$ in the
$L_u$-metric {\it for all} $1\le u\le\infty$. Thus
$$
\|f-P_Qf\|_{L_u(H_Q)}\approx
\inf_{P\in\Pc_{k-1}}\|f-P\|_{L_u(H_Q)}
$$
with constants depending only on $n,k$ and $\theta_S$. We
put $P_Qf:=0$ if $\diam Q>\delta_S$. We define an extension
$\tf$ by the formula
\bel{OE} \tf(x)=(\Eks f)(x):=\sum_{Q\in W_S}
\varphi_Q(x)P_Qf(x), \ \ \ \ \ x\in \RN\setminus S, \ee
and $\tf(x):=f(x),x\in S.$ Here $\{\varphi_Q:Q\in W_S\}$ is
a partition of unity subordinated to the Whitney
decomposition $W_S$.
\par This extension construction was first
used in \cite{S1,S2} to obtain a description of the
restrictions of Besov functions to regular sets. (In
Section 5 we present details of this approach and some main
facts related to the local approximation theory.) I am very
thankful to Yu.\ Brudnyi for the excellent suggestion that
the same construction might also yield a characterization
of the restriction of Triebel-Lizorkin functions to regular
sets via local approximations.
\par  We show that the extension operator $\tf=\Eks$
defined by \rf{OE} in some sense preserves local
approximation properties of functions, see Theorem
\ref{ExtLoc}. For example, Theorem \ref{Approx} states that
for every cube $Q=Q(x,r)$ such that $x\in S$ and
 $r\le\delta_S/4$ and every $1\le u\le\infty$
$$
\Ec_k(\tf;Q(x,r))_{L_u}
\le\gamma\Ec_k(f;Q(x,25r))_{L_u(S)}
$$
where $\gamma$ is a positive constant depending only on
$n,k,\theta_S$ and $\delta_S$.
\par In Section 4 we study extension properties of
certain generalized sharp maximal functions. These maximal
functions determine both the norm in the Sobolev space and
the norm in the Triebel-Lizorkin space. Given a vector of
parameters $\vv:=(s,k,q,u)$ where $0\le s\le k$,
$0<q\le\infty$, $1\le u\le\infty$ and a function $f\in
L_{u,\,loc}(S)$, we put
$$
\fs(x):= \left(\int_0^\infty
\left(\frac{\Ec_k(f;Q(x,t))_{L_u(S)}}{t^s}\right)^q\,
\frac{dt}{t}\right)^{\frac{1}{q}},\ \ \ \ \ x\in S,
$$
(with the corresponding modification if $q=\infty$).
\par Theorem \ref{MF} presents point-wise estimates of
$(\tf)^\sharp_{\vv,\RN}$ via the Hardy-Littlewood
ma\-xi\-mal function of $f$ and $\fs$. To formulate the
result given a function $h$ defined on $S$, we let $h^\cw$
denote its extension on all of $\RN$ by zero. We prove that
for every $1\le u\le\infty$, $0\le s<k$ (or $0\le s\le k$
if $q=\infty$) and every $x\in\RN$
\bel{Mfs} (\tf)^\sharp_{\vv,\RN}(x)\le
\gamma\{M((\fs)^\cw)(x)+M_u(f^\cw)(x)\} \ee
where $M_uf:=(M(|f|^u))^{\frac{1}{u}}$ and
$\gamma=\gamma(n,k,s,\theta_S,\delta_S)$. For instance,
$f^\sharp_\alpha=f^\sharp_{\vv,\RN}$ whenever
$\vv:=(\alpha,k,\infty,1)$, see \rf{SM}, so that by
\rf{Mfs} for every $0\le\alpha\le k$ we have
$$
(\tf)^{\sharp}_{\alpha}(x)\le
\gamma(M(f^{\sharp}_{\alpha,S})^\cw(x)+ Mf^\cw(x)),\ \ \ \
\ x\in\RN.
$$
\par Finally, applying the Hardy-Littlewood
maximal theorem to inequality \rf{Mfs} we obtain the
statements of Theorem \ref{EXT1} ($\vv:=(k,k,\infty,1)$)
and Theorem \ref{EXT2} ($\vv:=(s,k,q,1), 0<s<k$).
\par In turn the proof of Theorem \ref{EXT3} is based on
estimates of the modulus of continuity of the extension
$\tf$ via local approximations of $f$ on $S$. In Section 5
we prove that for every regular set $S$ and every function
$f\in L_p(S)$, $1\le p\le\infty,$ the modulus of continuity
of order $k$ of $\tf$ in $L_p$ satisfies the following
inequality
\bel{MS} \omega_k(\tf;t)_{L_p}\le \gamma\, t^k\left\{
\left( \int_t^{1}\left
(\frac{\|\Ec_k(f;Q(\cdot,\tau))_{L_p(S)}\|_{L_p(S)}}
{\tau^k}\right)^p \frac{d\tau}{\tau}\right)^\frac{1}{p}+
\|f\|_{L_p(S)}\right\}. \ee
Here $0<t\le 1$ and the constant $\gamma$ depends only on
$k,n,p,\theta_S$ and $\delta_S$.
\par Similar estimates for the quantity
$\|\Ec_k(\tf;Q(\cdot,t))_{L_u}\|_{L_p(\RN)}$ are given in
Theorem \ref{KPS}. Using these estimates, description
\rf{ChB} and the Hardy inequality we obtain the result of
Theorem \ref{EXT3}.
\begin{remark} Observe that  the operator
$\Eks$ defined by \rf{OE} provides a ``universal" {\it
linear continuous extension operator} from $\WKP|_S$,
$\FS|_S$ and $\BS|_S$ into corresponding spaces on $\RN$.
(The ``universality"  means that, for all sufficiently
large $k$, the operator $\Eks$ depends only on the regular
set $S$ and is independent of the indices of the spaces).
This allows us to complement the statements of Theorems
\ref{EXT1}, \ref{EXT2} and \ref{EXT3} with the following
assertion:
\par {\it There exists a linear extension
operator mapping functions on $S$ to functions on $\RN$
which is continuous from $\Ac|_S$ into $\Ac$ whenever $\Ac$
is any one of the following spaces:
\par (i) $\WKP$ for $1<p\le\infty$,
\par (ii) $\FS$ for  $s>0,$ $1<p<\infty$,
and $1\le q\le\infty$,
\par (iii) $\BS$ for  $s>0, 1\le p\le\infty,$
$0< q\le\infty$.
\par The norm of this operator is bounded by a constant
which depends only on $n,\theta_S,\delta_S$ and the
parameters of the space $\Ac$.}
\par For the Triebel-Lizorkin spaces
the existence of a linear continuous extension ope\-ra\-tor
from $\FS|_S$ into $\FS$ where $S$ is an arbitrary regular
set, has been proved by Rychkov \cite{R1}. For the scale of
the Besov spaces this follows from a result of Shvartsman
\cite{S1}, see also \cite{S2}.
\end{remark}
\begin{remark} As we have noted above, our goal is to give
a {\it constructive} intrinsic characterization of the
restrictions of smooth functions to regular sets. The
descriptions of trace spaces given in Theorems \ref{EXT1},
\ref{EXT2} and \ref{EXT3} are not quite constructive,
because they use the notion of the best local polynomial
approximation of a function on a regular set, the quantity
$\Ec_k(f;Q)_{L_u(S)}$, see \rf{ES}.
\par However, following an idea of Yu. Brudnyi
\cite{Br4}, we can readily eliminate this element of
nonconstructi\-vi\-ty. In fact, using Proposition
\ref{LinPr} and the regularity of $S$, we immediately
obtain the equivalence
$$
\Ec_k(f;Q(x,r))_{L_u(S)}\approx |Q|^{-\frac{1}{u}}
\|f-\PR_{k,Q\cap S}(f)\|_{L_u(Q\cap S)}, \ \ \ \ x\in S,
 \ \ r\le 1.
$$
Here $1\le u\le \infty$ and $\PR_{k,Q\cap
S}(f)\in\Pc_{k-1}$ denotes the polynomial of the best
approximation of $f$ on $Q\cap S$ {\it in $L_2$-norm}. Of
course there are many constructive formulas for calculation
of this polynomial, see e.g. \rf{PRL2}. For instance, for
$k=1$, one can take $\PR_{1,Q\cap S}(f)$ to be the average
of $f$ over $Q\cap S$ so that in this case
$$
\Ec_1(f;Q(x,r))_{L_u(S)}\approx \left(\frac{1}{|Q|}
\int_{Q\cap S}
 \left|f-\frac{1}{|Q\cap S|}\int_{Q\cap S}fdy
 \right|^u\,dx\right)
 ^{\frac{1}{u}},\ \ \ \ \ x\in S, \ \ r\le 1.
$$
\end{remark}
\SECT{2. The Whitney covering and a family of associated
``quasi-cubes".}{2}
\indent
\par Our notation is fairly standard.
Throughout the paper $C,C_1,C_2,...$ or $\gamma,
\gamma_1,\gamma_2,...$ will be generic positive constants
which depend only on $n,\theta_S,\delta_S$ and indexes of
spaces ($s,p,q,k,$ etc.). These constants can change even
in a single string of estimates. The dependence of a
constant on certain parameters is expressed, for example,
by the notation $\gamma=\gamma(n,k,p)$. We write $A\approx
B$ if there is a constant $C\ge 1$ such that $A/C\le B\le
CA$.
\par It will be convenient for us to measure distances in
$\RN$ in the uniform norm
$$
\|x\|_\infty:=\max\{|x_i|:~i=1,...,n\}, \ \ \
x=(x_1,...,x_n)\in\RN.
$$
Thus every cube
$$
Q=Q(x,r):=\{y\in\RN:\|y-x\|_\infty\le r\}
$$
is a ``ball" in $\|\cdot\|_\infty$-norm  of ``radius" $r$
centered at $x$. We let $x_Q:=x$ denote center of $Q$ and
$r_Q:=r$ its ``radius". Given a constant $\lambda> 0$, we
let $\lambda Q$ denote the cube $Q(x,\lambda r)$. By $Q^*$
we denote the cube $Q^*:=\frac{9}{8}Q$.
\par As usual given subsets $A,B\subset \RN$,  we put
$
\diam A:=\sup\{\|a-a'\|_\infty:~a,a'\in A\}
$
and
$$
\dist(A,B):=\inf\{\|a-b\|_\infty:~a\in A, b\in B\}.
$$
We also set $\dist(x,A):=\dist(\{x\},A)$ for $x\in \RN$. By
$\cl(A)$ we denote the closure of $A$ in $\RN$. Finally,
$\chi_A$ denotes the characteristic function of $A$; we put
$\chi_A\equiv 0$ if $A=\emptyset$.
\par The following property of regular sets is well known
(see, e.g. \cite{S3}).
\begin{lemma}
$|\cl(S)\setminus S|=0$ provided $S$ is a regular subset of
$\RN$.
\end{lemma}
\par In the remaining part of the paper we will assume
that $S$ is a {\it closed} regular subset of $\RN$. Since
now $\RN\setminus S$ is an open set, it admits a Whitney
decomposition $W_S$, see, e.g. Stein \cite{St}. We recall
the main properties of $W_S$.
\begin{theorem}\label{Wcov}
$W_S=\{Q_k\}$ is a countable family of cubes  such that
\par (i). $\RN\setminus S=\cup\{Q:Q\in W_S\}$;
\par (ii). For every cube $Q\in W_S$
$$
\diam Q\le \dist(Q,S)\le 4\diam Q;
$$
\par (iii). Every point of $\RN\setminus S$ is covered by
at most $N=N(n)$ cubes from $W_S$.
\end{theorem}
\par We also need certain additional properties of Whitney's
cubes which we present in the next lemma. These properties
readily follow from (i)-(iii).
\begin{lemma}\label{Wadd}
(1). If $Q,K\in W_S$ and $Q^*\cap K^*\ne\emptyset$, then
$$
\frac{1}{4}\diam Q\le \diam K\le 4\diam Q.
$$
(Recall that $Q^*:=\frac{9}{8}Q$.)
\par (2). For every cube $K\in W_S$ there are at most
$N=N(n)$ cubes from the family
$
W^*_S:=\{Q^*:Q\in W_S\}
$
which intersect $K^*$.
\end{lemma}
\par Observe that the family of cubes $W_S$ constructed in
\cite{St} satisfies conditions of Theorems \ref{Wcov} and
Lemma \ref{Wadd} with respect to the Euclidean norm rather
than the uniform one. However, a simple modification of
this construction provides a family of Whitney's cubes
which have the same properties with respect to the uniform
norm.
\par Let us formulate the main result of the section.
\begin{theorem}\label{HB} Let $S$ be a regular subset of
$\RN$. There is a family of Borel sets $\Hc_S=\{H_Q:~Q\in
W_S\}$ such that:
\par (i). $H_Q\subset (10Q)\cap S,\ \ \ \ Q\in W_S$;
\par (ii). $|\,Q|\le \gamma_1|\,H_Q|$ whenever
$Q\in W_S$ and $\diam Q\le\delta_S$;
\par (iii). $\sum_{Q\in W_S}\limits\chi_{H_Q}\le \gamma_2$.
\par Here  $\gamma_1,\gamma_2$ are positive
constants depending only on $n$ and $\theta_S$.
\end{theorem}
{\bf Proof.} Let $K=Q(x_K,r_K)\in W_S$ and let $a_K\in S$
be a point nearest to $x_K$ on $S$. Then by property (ii)
of Theorem \ref{Wcov}
\bel{inBY}
Q(a_K,r_K)\subset 10K.
\ee
\par Given $\ve, 0<\ve\le 1,$, we denote
$K_{\ve}:=Q(a_K,\ve r_K)$. Let $Q=Q(x_Q,r_Q)$ be a cube
from $W_S$ with $\diam Q\le \delta_S$. Set
\bel{2.13'}
\Ac_Q:=\{K=Q(x_K,r_K)\in W_S:~K_\ve\cap
Q_\ve\ne\emptyset,~r_K\le \ve r_Q\}.
\ee
(Recall that $Q_{\ve}:=Q(a_Q,\ve r_Q)$.) We define a
``quasi-cube" $H_Q$ by letting
\bel{dAB}
H_Q:=(Q_\ve\cap S)\setminus (\cup\{K_\ve:~K\in
\Ac_Q\}).
\ee
If $\diam Q>\delta_S$ we put $H_Q:=\emptyset$.
\par Prove that for some $\ve:=\ve(n,\theta_S)\le 1$
small enough the family of subsets $\Hc_S$ satisfies
conditions (i)-(iii). By \rf{dAB} and \rf{inBY}
\bel{QE1}
H_Q\subset Q_\ve:=Q(a_Q,\ve r_Q)
 \subset Q(a_Q,r_Q)\subset 10Q.
\ee
In addition, by \rf{dAB} $H_Q\subset S$ so that $H_Q\subset
(10Q)\cap S$ proving property (i).
\par Let us prove (ii). Suppose that $Q=Q(x_Q,r_Q)\in W_S$
and $\diam Q\le\delta_S$. If $K\in \Ac_Q$, then by
\rf{2.13'} $K_\ve\cap Q_\ve\ne \emptyset$ and $r_K\le\ve
r_Q$. Hence
$$
r_{K_\ve}(=\ve r_K)\le\ve r_{Q_\ve}(=\ve^2r_Q) \le
r_{Q_\ve}
$$
so that $a_K\in 2Q_\ve$. Since  $Q(a_K,r_K)\subset 10K$,
see \rf{inBY}, $K\subset Q(a_K,10r_K)$ as well. In
addition, $r_K\le\ve r_Q=r_{Q_\ve}$ which implies $K\subset
12Q_\ve$. Thus
\bel{2.**}
\Uc_Q:=\cup\{K:~K\in \Ac_Q\}\subset 12Q_\ve.
\ee
By property (iii) of Theorem \ref{Wcov}
$$
\sum_{K\in \Ac_Q}\limits\chi_K\le \sum_{K\in
W_S}\limits\chi_K\le N=N(n),
$$
so that by \rf{2.**}
$$
\sum_{K\in \Ac_Q}|\,K|= \int\limits_{\Uc_Q} \sum_{K\in
\Ac_Q}\chi_K\, dx \le\int\limits_{12Q_\ve}
Ndx=N|\,12Q_\ve|= N12^n|\,Q_\ve|= C_1|\,Q_\ve|.
$$
On the other hand, for every $K=Q(x_K,r_K)\in \Ac_Q$ we
have
$$
|\,K_\ve|=|\,Q(a_K,\ve r_K)|=\ve^n|\,Q(a_K,r_K)|=
\ve^n|\,K|.
$$
Hence
$$
|\,\cup\{K_\ve:~K\in \Ac_Q\}| \le \sum_{K\in
\Ac_Q}|\,K_\ve| =\ve^n\sum_{K\in \Ac_Q}|\,K| \le
C_1\ve^n|\,Q_\ve|.
$$
\par Since $S$ is regular and
$\diam Q_\ve=\ve \diam Q\le \delta_S$, $|\,Q_\ve\cap S|\ge
\theta_S^{-1} |\,Q_\ve|$ so that
\be |\,H_Q|&=&|\,(Q_\ve\cap S|
\setminus (\cup\{K_\ve:~K\in \Ac_Q\})|\nn\\
&\ge& |\,Q_\ve\cap S|- |\,\cup\{K_\ve:~K\in \Ac_Q\}|\ge
(\theta_S^{-1}-C_1\ve^n)|\,Q_\ve|.\nn \ee
Clearly,
$ |\,Q_\ve|=|\,Q(a_Q,\ve r_Q)|=\ve^nr_Q^n=\ve^n|\,Q|$
so that
$$
|\,H_Q|\ge (\theta_S^{-1}-C_1\ve^n)\ve^n|\,Q|.
$$
We define $\ve$ by setting
$\ve:=(2C_1\theta_S)^{-\frac{1}{n}}$. Then the inequality
$|\,Q|\le\gamma_1|\,H_Q|$ holds with
$\gamma_1:=4C_1\theta_S^2$ proving property (ii) of the
theorem.
\par Let us prove (iii). Let
$Q=Q(x_Q,r_Q),Q'=Q(x_{Q'},r_{Q'})\in W_S$ be Whitney's
cubes such that $\diam Q$, $\diam Q'\le\delta_S$ and $H_Q\cap
H_{Q'}\ne \emptyset$. Since $H_Q\subset Q_\ve,
H_{Q'}\subset Q'_\ve$, we have $Q_\ve\cap Q'_\ve\ne
\emptyset$.
\par On the other hand,
$Q\notin \Ac_{Q'}$ and $Q'\notin \Ac_{Q}$, otherwise by
\rf{2.13'} and \rf{dAB} $H_Q\cap H_{Q'}=\emptyset$. Since
$Q_\ve\cap Q'_\ve\ne \emptyset$, by definition \rf{2.13'}
$r_Q>\ve r_{Q'}$ and $r_{Q'}>\ve r_Q$ so that $r_Q\approx
r_{Q'}$. By \rf{QE1} $Q_\ve\subset 10Q $ and similarly
$Q'_\ve\subset 10Q'$. But $Q_\ve\cap Q'_\ve\ne \emptyset$
so that $10Q\cap 10Q'\ne\emptyset$ as well. Since
$r_Q\approx r_{Q'}$, this imply $Q'\subset C_2Q$ for some
constant $C_2=C_2(\ve)=C_2(n,\theta_S)$. Observe also that
$|\,Q|\approx |\,Q'|$.
\par We denote
$$
\Tc_Q:=\{Q'\in W_S:~H_Q\cap H_{Q'} \ne\emptyset, ~\diam
Q'\le\delta_S\}
$$
and $\Vc_Q:=\cup\{Q':~Q'\in\Tc_Q\}$. Thus we have proved
that $\Vc_Q\subset C_2Q$ and $|\,Q'|\approx |\,Q|$ for
every $Q'\in\Tc_Q$.
\par Let $M_Q:=\card\Tc_Q$ be the cardinality of $\Tc_Q$.
Clearly, to prove (iii) it suffices to show that
$M_Q\le\gamma_2$. We have
$$
M_Q|\,Q|\le C\sum_{Q'\in\Tc_Q}|\,Q'|= C\int_{\Vc_Q}\limits
\sum_{Q'\in\Tc_Q}\chi_{Q'}\,dx \le C\int_{C_2Q}\limits
\sum_{Q'\in\Tc_Q}\chi_{Q'}\,dx.
$$
By the property (iii) of Theorem \ref{Wcov}
$$
\sum\{\chi_{Q'}:Q'\in\Tc_Q\}\le \sum\{\chi_{Q'}:Q'\in
W_S\}\le N=N(n)
$$
so that
$$
M_Q|\,Q|\le C\int_{C_2Q}\limits Ndx=C N|\,C_2Q| =C N
C_2^n|\,Q|
$$
proving the required inequality $M_Q\le\gamma_2$. \bx
\SECT{3. Local approximation properties of the extension
operator.}{3}
\indent
\par In this section we present estimates of local
polynomial approximations of the extension $\tf$, see
\rf{OE}, via corresponding local approximation of a
function $f$ defined on a regular subset $S\subset\RN$. We
start by presenting two lemmas about properties of
polynomials on subsets of $\RN$.
\begin{proposition}\label{Pol}
(Brudnyi and Ganzburg \cite{BrG}) Let $A$ be a measurable
subset of a cube $Q$, $|A|>0$, $1\le u_1,u_2\le \infty$ and
$P\in\Pc_k$. Then
$$
|Q|^{-\frac{1}{u_1}}\|P\|_{L_{u_1}(Q)}\le\gamma
|A|^{-\frac{1}{u_2}}\|P\|_{L_{u_2}(A)}
$$
where $\gamma$ is a positive constant depending only on
$n,k$ and the ratio $|Q|/|A|$.
\end{proposition}
\par The proposition implies two corollaries.
\begin{corollary}\label{CorA}
For every subset $A$ of a cube $Q$, $|A|>0$, every $1\le
u_1,u_2\le \infty$ and every polynomial $P\in\Pc_k$
\bel{DifM} |A|^{-\frac{1}{u_1}}\|P\|_{L_{u_1}(A)}\le\gamma
|A|^{-\frac{1}{u_2}}\|P\|_{L_{u_2}(A)} \ee
where $\gamma$ depends only on $n,k$ and $|Q|/|A|$.
\end{corollary}
\begin{corollary}\label{CorB}
Let $A_i\subset Q_i$, $|A_i|>0$, $i=1,2$. Suppose that
$(\lambda_1Q_1)\cap(\lambda_1Q_2)\ne\emptyset$ and
$\lambda_2^{-1}r_{Q_1}\le r_{Q_2}\le\lambda_2r_{Q_1}$ where
$\lambda_1,\lambda_2$ are some positive constants. Then for
every $1\le u\le \infty$ and every polynomial $P\in\Pc_k$
$$
\|P\|_{L_u(A_1)}\le\gamma \|P\|_{L_u(A_2)}
$$
where $\gamma$ depends only on $n,k,\lambda_i$ and
$|Q_i|/|A_i|$, $i=1,2$.
\end{corollary}
\par Given a function $f\in L_{u,\,loc}(\RN)$,
$1\le u\le\infty$, and a measurable subset $A\subset\RN$,
we let $E_k(f;A)_{L_u}$ denote the {\it local best
approximation} of order $k$ of $f$ on $A$ in $L_u$-norm,
see Brudnyi \cite{Br1},
\bel{EQ}
E_k(f;A)_{L_u}:=\inf_{P\in\Pc_{k-1}}
\|f-P\|_{L_u(A)}.
\ee
Thus
$$
\Ec_k(f;Q)_{L_u(S)}=|Q|^{-\frac{1}{u}}E_k(f;Q\cap S)_{L_u}
$$
see \rf{ES}. We note a simple property of
$\Ec_k(f;\cdot)_{L_u(S)}$ as a cube function: for every two
cubes $Q_1\subset Q_2$
\bel{ESUB} \Ec_k(f;Q_1)_{L_u(S)}\le
\left(\frac{|Q_2|}{|Q_1|}\right)
^{\frac{1}{u}}\Ec_k(f;Q_2)_{L_u(S)}. \ee
\begin{proposition}\label{LinPr}
(Brudnyi \cite{Br4}) Let $A$ be a subset of a cube $Q$,
$|A|>0$. Then there is a linear operator
$\PR_{k,A}:L_1(A)\to \Pc_{k-1}$ such that for every $1\le
u\le \infty$ and every $f\in L_u(A)$
$$
\|f-\PR_{k,A}(f)\|_{L_u(A)}\le\gamma E_k(f;A)_{L_u}.
$$
Here $\gamma=\gamma(n,k,\frac{|Q|}{|A|})$.
\end{proposition}
{\bf Proof.} Recall the construction of $\PR_{k,A}$ given
in \cite{Br4}. We let $\{P_\beta:|\beta|\le k-1\}$ denote
an orthonormal basis in the linear space $\Pc_{k-1}$ with
respect to the inner product $\ip{f,g}=\int_A
f(x)g(x)\,dx$. We put
\bel{PRL2}
\PR_{k,A}(f):=\sum_{|\beta|\le k-1}\left(\int_A
P_\beta(x) f(x)\,dx\right)P_\beta.
\ee
Clearly, $\PR_{k,A}:L_1(A)\to \Pc_{k-1}$ is a projector
(i.e., $\PR_{k,A}^2=\PR_{k,A}$). Estimate its operator norm
in $L_u(A)$. For every $f\in L_u(A)$ we have
$$
\|\PR_{k,A}(f)\|_{L_u(A)}\le \sum_{|\beta|\le
k-1}\left|\int_A P_\beta(x)
f(x)\,dx\right|\|P_\beta\|_{L_u(A)}
$$
so that by the H\"{o}lder inequality
$$
\|\PR_{k,A}(f)\|_{L_u(A)}\le\left(\sum_{|\beta|\le k-1}
\|P_\beta\|_{L_{u}(A)}\|P_\beta\|_{L_{u^*}(A)}\right)
\|f\|_{L_u(A)}
$$
where $1/u+1/u^*=1$. But by \rf{DifM}
$$
\|P_\beta\|_{L_u(A)}\|P_\beta\|_{L_{u^*}(A)}\le \gamma^2
(|A|^{\frac{1}{u}-\frac{1}{2}}\|P_\beta\|_{L_2(A)})
(|A|^{\frac{1}{u^*}-\frac{1}{2}}\|P_\beta\|_{L_2(A)})
=\gamma^2
$$
proving that
$\|\PR_{k,A}(f)\|_{L_u(A)}\le\gamma_1\|f\|_{L_u(A)}$ with
$\gamma_1=\card\{\beta:|\beta|\le k-1\}\,\gamma^2$ (recall
that $\|P_\beta\|_{L_2(A)}=1$ for every $\beta$). The last
inequality in the standard way implies
$$
\|f-\PR_{k,A}(f)\|_{L_u(A)}\le (1+\|\PR_{k,A}\|)
E_k(f;A)_{L_u}\le (1+\gamma_1) E_k(f;A)_{L_u}.\BX
$$
\par Proposition \ref{LinPr} and Theorem \ref{HB}
immediately imply the following
\begin{corollary}\label{PrQ}
Let $S$ be a regular set and let $Q\in W_S$ be a cube with
$\diam Q\le\delta_S$. There is a linear continuous operator
$P_Q:L_1(H_Q)\to \Pc_{k-1}$ such that for every function
$f\in L_{u,\,loc}(S)$, $1\le u\le\infty,$
$$
\|f-P_Qf\|_{L_u(H_Q)}\le \gamma E_k(f;H_Q)_{L_u}.
$$
Here $\gamma=\gamma(n,k,\theta_S)$.
\end{corollary}
\par We put
\bel{PO} P_Qf=0,\ \ \ \ {\rm if}\ \ \ \ \diam Q>\delta_S.
\ee
Now the map $Q\to P_Q(f)$ is defined on all of the family
$W_S$.  This map gives rise a linear extension operator
defined by the formula
\bel{ExtOp}
\Eks f(x):=\left \{
\begin{array}{ll}
f(x),& x\in S,\\\\
\sum\limits_{Q\in W_S} \varphi_Q(x)P_Qf(x),& x\in
\RN\setminus S.
\end{array}
\right.
\ee
\par Here $\Phi_S:=\{\varphi_Q:Q\in W_S\}$
is a smooth partition of unity subordinated to the Whitney
decomposition $W_S$, see, e.g. \cite{St}. We recall that
$\Phi_S$ is a family of functions defined on $\RN$ which
have the following properties:
\par (a). $0\le\varphi_Q\le 1$ for every $Q\in W_S$;
\par (b). $\supp \varphi_Q\subset Q^*(:=\frac{9}{8}Q),$
$Q\in W_S$;
\par (c). $\sum\{\varphi_Q(x):~Q\in W_S\}=1$ for every
$x\in\RN\setminus S$;
\par (d). for every multiindex $\beta, |\beta|\le k$ and
every cube $Q\in W_S$
$$
|D^\beta\varphi_Q(x)| \le C (\diam Q)^{-|\beta|}, \ \ \ \
x\in\RN,
$$
where $C$ is a constant depending only on $n$ and $k$.
\par We turn to estimates of local approximations of the
extension operator
$$
\tf:=\Eks f.
$$
To formulate the main result of the section, Theorem
\ref{ExtLoc}, given $x\in\RN$ and $t>0$ we let $a_x$ denote
a point nearest to $x$ on $S$ (in the uniform metric). Thus
$\|x-a_x\|_\infty=\dist(x,S)$. We put
\bel{RXT} r^{(x,t)}:=50\max(80t,\dist(x,S)) \ee
and
\bel{KXT} K^{(x,t)}:=Q(a_x,r^{(x,t)}). \ee
\begin{theorem}\label{ExtLoc}
Let $S$ be a regular subset of $\RN$ and let $f\in
L_{u,\,loc}(S)$, $1\le u\le \infty$. Then for every
$x\in\RN$ and $t>0$
$$
\Ec_k(\tf;Q(x,t))_{L_u}\le C\frac{t^k}{t^k+\dist(x,S)^k}
\left \{
\begin{array}{ll}
\Ec_k(f;K^{(x,t)})_{L_u(S)},& r^{(x,t)}\le\delta_S,\\\\
\Ec_0(f;K^{(x,t)})_{L_u(S)},& r^{(x,t)}>\delta_S.
\end{array}
\right.
$$
Here $\gamma=\gamma(n,k,\theta_S,\delta_S)$.
\end{theorem}
\par We recall that $\Pc_{-1}:=\{0\}$ so that by definition
\rf{ES}
\bel{K0}
\Ec_0(f;K^{(x,t)})_{L_u(S)}:=
|K^{(x,t)}|^{-\frac{1}{u}}\|f\|_{L_u(K^{(x,t)}\cap S)}.
\ee
\par We will prove the theorem for the case
$1\le u<\infty$; corresponding changes for $u=\infty$ are
obvious.
\par The proof is based on a series of auxiliary lemmas.
To formulate the first of them given a cube $K\subset\RN$,
we define two families of Whitney's cubes:
$$
\Qc_1(K):=\{Q\in W_S:~Q\cap K\ne\emp\}
$$
and
\bel{AA} \Qc_2(K):=\{Q\in W_S:\exists\, Q'\in\Qc_1(K)\ \ \
{\rm such~that~~}Q'\cap Q^*\ne\emp\}. \ee
\begin{lemma}\label{Cubes} Let $K$ be a cube centered
in $S$. Then for every $Q\in\Qc_2(K)$  we have $\diam Q\le
2\diam K$ and $\|x_K-x_Q\|_\infty \le \frac{5}{2}\diam K.$
\end{lemma}
{\bf Proof.} Since $Q\in\Qc_2(K)$, there is a cube
$Q'\in\Qc_1(K)$ such that $Q'\cap Q^*\ne\emp$. Let $a\in
Q'\cap K$. Since $x_K\in S$, by property (ii) of Theorem
\ref{Wcov}
$$
\diam Q'\le \dist(Q',S)\le\dist(a,S)\le\|a-x_K\|_\infty\le
\frac{1}{2}\diam K
$$
so that $\diam Q'\le \frac{1}{2}\diam K.$
But $\diam Q\le 4\diam Q'$, see Lemma \ref{Wadd}, (1),
proving that $\diam Q\le 2\diam K.$
\par Recall that $Q'\cap K\ne\emp$, $Q'\cap Q^*\ne\emp$ and
$\diam Q'\le \frac{1}{2}\diam K.$ It remains to make use of
the triangle inequality and the required inequality
$\|x_K-x_Q\|_\infty \le (5/2)\diam K$ follows. \bx
\begin{lemma}\label{EONE}
Let $S$ be a regular subset of $\RN$ and let $f\in
L_{u,\,loc}(S)$, $1\le u\le\infty$. Then for every cube $K$
and every polynomial $P_0\in\Pc_{k-1}$
$$
\|\tf-P_0\|^u_{L_u(K\setminus S)}\le
C\sum_{Q\in\Qc_2(K)} \|P_{Q}-P_0\|^u_{L_u(Q)}.
$$
\end{lemma}
{\bf Proof.} Clearly, $K\setminus S\subset
\cup\{Q:~Q\in\Qc_1(K)\}$ so that
$$
\|\tf-P_0\|^u_{L_u(K\setminus S)}\le \sum_{Q\in\Qc_1(K)}
\|\tf-P_0\|^u_{L_u(Q)}.
$$
\par Let $Q\in\Qc_1(K)$. By $V(Q)$ we denote a family of
cubes
$$
V(Q):=\{Q'\in W_S: (Q')^*\cap Q\ne\emp\}.
$$
Clearly, by property (2) of Lemma \ref{Wadd} $M_Q:=\card
V(Q)\le N(n)$. Properties (a)-(c) of the partition of unity
and formula \rf{ExtOp} imply
\be \|\tf-P_0\|^u_{L_u(Q)}&\le& \|\sum_{Q'\in W_S}
\varphi_{Q'}(P_{Q'}-P_0)\|^u_{L_u(Q)}\nn\\&=& \|\sum_{Q'\in
V(Q)} \varphi_{Q'}(P_{Q'}-P_0)\|^u_{L_u(Q)} \le
M_Q^{u-1}\sum_{Q'\in V(Q)} \|P_{Q'}-P_0\|^u_{L_u(Q)} \nn
\ee
so that
$$
\|\tf-P_0\|^u_{L_u(Q)}\le C\sum_{Q'\in V(Q)}
\|P_{Q'}-P_0\|^u_{L_u(Q)}.
$$
Since $(Q')^*\cap Q\ne\emp$ for every $Q'\in V(Q)$, by
Lemma \ref{Wadd}, (1), $\diam Q'\approx \diam Q$. Hence by
Corollary \ref{CorB}
$$
\|P_{Q'}-P_0\|_{L_u(Q)}\approx \|P_{Q'}-P_0\|_{L_u(Q')}
$$
so that
$$
\|\tf-P_0\|^u_{L_u(Q)}\le C\sum_{Q'\in V(Q)}
\|P_{Q'}-P_0\|^u_{L_u(Q')}.
$$
This implies
$$
\|\tf-P_0\|^u_{L_u(K\setminus S)}\le C\sum_{Q\in\Qc_1(K)}
\sum_{Q'\in V(Q)}\|P_{Q'}-P_0\|^u_{L_u(Q')}.
$$
Clearly, every cube $Q'$ on the right-hand side of this
inequality belongs to $\Qc_2(K)$, see definition \rf{AA}.
Moreover, by Lemma \ref{Wadd}, (2), for every such a cube
$Q'$ there are at most $N(n)$ cubes $Q\in W_S$ such that
 $V(Q)\ni Q'$. Hence
\be \|\tf-P_0\|^u_{L_u(K\setminus S)}&\le&
C\sum_{Q'\in\Qc_2(K)} \card\{Q: V(Q)\ni Q'\}
\|P_{Q'}-P_0\|^u_{L_u(Q')}\nn\\
&\le& C N(n)\sum_{Q\in\Qc_2(K)}
\|P_{Q}-P_0\|^u_{L_u(Q)}.\BX\nn
\ee
\par Given a cube $K\subset\RN$, define a family of cubes
\bel{AA3}
\Qc_3(K):=\{Q\in \Qc_2(K):~\diam Q\le \delta_S\}.
\ee
\begin{lemma}\label{ETWO}
Let $S$ be a regular subset of $\RN$ and let $f\in
L_{u,\,loc}(S)$, $1\le u\le\infty$. Then for every cube $K$
centered in $S$ and every polynomial $P_0\in\Pc_{k-1}$
$$
\sum_{Q\in\Qc_3(K)} \|P_{Q}-P_0\|^u_{L_u(Q)} \le
C\|f-P_0\|^u_{L_u((25K)\cap S)}.
$$
\end{lemma}
{\bf Proof.} For each $Q\in\Qc_3(K)$ by properties (i),(ii)
of Theorem \ref{HB} and by Corollary \ref{CorB} we have
$$
\|P_{Q}-P_0\|_{L_u(Q)}\le C \|P_{Q}-P_0\|_{L_u(H_Q)}.
$$
By Corollary \ref{PrQ}
$$
\|P_{Q}-P_0\|_{L_u(H_Q)}\le \|P_{Q}-f\|_{L_u(H_Q)}
+\|f-P_0\|_{L_u(H_Q)}\le \gamma
E_k(f;H_Q)_{L_u}+\|f-P_0\|_{L_u(H_Q)}
$$
where $\gamma=\gamma(n,k,\theta_S)$. Since
$E_k(f;H_Q)_{L_u}\le\|f-P_0\|_{L_u(H_Q)},$
see definition \rf{EQ}, we have
$$
\|P_{Q}-P_0\|_{L_u(Q)}\le C \|f-P_0\|_{L_u(H_Q)}.
$$
Put $B:=\cup\{H_Q:Q\in\Qc_3(K)\}$ and $\eta:=\sum
\{\chi_{H_Q}:Q\in\Qc_3(K)\}$. Then the last inequality
imply
$$
\sum_{Q\in\Qc_3(K)} \|P_{Q}-P_0\|^u_{L_u(Q)}\le C
\sum_{Q\in\Qc_3(K)}\|f-P_0\|^u_{L_u(H_Q)}
=C\|\eta\,(f-P_0)\|^u_{L_u(B)}.
$$
But by property (iii) of Theorem \ref{HB}~
$\eta\le\sum\{\chi_{H_Q}:Q\in W_S\} \le \gamma(n,\theta_S)$
so that
$$
\sum_{Q\in\Qc_3(K)} \|P_{Q}-P_0\|^u_{L_u(Q)}\le
C\|f-P_0\|^u_{L_u(B)}.
$$
By Lemma \ref{Cubes} for each $Q\in\Qc_3\subset\Qc_2$ we
have $\|x_K-x_Q\|_\infty \le (5/2)\diam K$ and $\diam Q\le
2\diam K$. Moreover, by property (i) of Theorem \ref{HB},
$H_Q\subset (10Q)\cap S$. Hence
$$
H_Q\subset 10 Q \subset (10\cdot 2+5)K=25K
$$
proving that $B\subset (25 K)\cap S$.\bx
\begin{proposition}
Let $f\in L_{u,\,loc}(S)$, $1\le u\le\infty$, where $S$ is
a regular set. Then for every cube $K$ with $\diam K\le
\delta_S/2$ centered in $S$ and every polynomial
$P_0\in\Pc_{k-1}$
$$
\|\tf-P_0\|_{L_u(K)}\le C\|f-P_0\|_{L_u((25K)\cap S)}.
$$
\end{proposition}
{\bf Proof.} By Lemma \ref{EONE}
$$
\|\tf-P_0\|^u_{L_u(K\setminus S)}\le C\sum_{Q\in\Qc_2(K)}
\|P_{Q}-P_0\|^u_{L_u(Q)}.
$$
Since $\diam K\le\delta_S/2$, by Lemma \ref{Cubes} for
every $Q\in\Qc_2(K)$  we have $\diam Q\le 2\diam K$ so that
$\diam Q\le \delta_S$. Hence  $\Qc_2(K)=\Qc_3(K)$, see
definition \rf{AA3}.
\par Therefore by Lemma \ref{ETWO}
$$
\sum_{Q\in\Qc_2(K)} \|P_{Q}-P_0\|^u_{L_u(Q)}=
\sum_{Q\in\Qc_3(K)} \|P_{Q}-P_0\|^u_{L_u(Q)}\le
C\|f-P_0\|^u_{L_u((25K)\cap S)}
$$
so that
$$
\|\tf-P_0\|^u_{L_u(K\setminus S)}\le
C\|f-P_0\|^u_{L_u((25K)\cap S)}.
$$
Finally,
$$
\|\tf-P_0\|^u_{L_u(K)}= \|\tf-P_0\|^u_{L_u(K\cap
S)}+\|\tf-P_0\|^u_{L_u(K\setminus S)}\le
(C+1)\|f-P_0\|^u_{L_u((25K)\cap S)}
$$
proving the lemma.\bx
\par Let us put $P_0\in\Pc_{k-1}$ to be a polynomial of
the best approximation of $f$ on $(25K)\cap S$ in
$L_u$-norm. Then the above proposition implies the
following inequality
$$
E_k(\tf;K)_{L_u}\le C E_k(f;(25K)\cap S)_{L_u}.
$$
Since $|K|\approx |25K|$, we obtain the next
\begin{theorem}\label{Approx} Let $S$ be a regular set and
let $f\in L_{u,\,loc}(S)$, $1\le u\le\infty$. Then for
every cube $K$ with $\diam K\le\delta_S/2$ centered in $S$
$$
\Ec_k(\tf;K)_{L_u}\le C \Ec_k(f;25K)_{L_u(S)}.
$$
\end{theorem}
\par Let us estimate the $L_u$-norm of the extension $\tf$.
\begin{proposition}\label{ELU}
Let $f\in L_{u,\,loc}(S)$, $1\le u\le\infty$, where $S$ is
regular. Then for every cube $K$ centered in $S$
$$
\|\tf\|_{L_u(K)}\le C\|f\|_{L_u((25K)\cap S)}.
$$
\end{proposition}
{\bf Proof.} By Lemma \ref{EONE} with $P_0:=0$ we have
$$
\|\tf\|^u_{L_u(K\setminus S)}\le C\sum_{Q\in\Qc_2(K)}
\|P_{Q}\|^u_{L_u(Q)}.
$$
Recall that $P_Q:=0$ if $\diam Q>\delta_S$, see \rf{PO}, so
that by definition \rf{AA3}
$$
\|\tf\|^u_{L_u(K\setminus S)}\le C\sum_{Q\in\Qc_3(K)}
\|P_{Q}\|^u_{L_u(Q)}.
$$
Now by Lemma \ref{ETWO} (with $P_0=0$) we obtain
$$
\sum_{Q\in\Qc_3(K)} \|P_{Q}\|^u_{L_u(Q)} \le
C\|f\|^u_{L_u((25K)\cap S)}
$$
so that
$\|\tf\|_{L_u(K\setminus S)}\le C\|f\|_{L_u((25K)\cap S)}.$
Finally,
$$
\|\tf\|^u_{L_u(K)}=\|\tf\|^u_{L_u(K\cap S)}+
\|\tf\|^u_{L_u(K\setminus S)} \le
(C+1)\|f\|^u_{L_u((25K)\cap S)}.\BX
$$
\par We turn to estimates of local approximations of $\tf$
on cubes which are located rather far from the set $S$. In
the remaining part of the section we will assume that a
cube $K=Q(x_K,r_K)$ satisfies the inequality
\bel{KA} \diam K\le\dist(x_K,S)/40. \ee
We let $Q_K\in W_K$ denote a Whitney's cube which contains
center of $K$, the point $x_K$.
\begin{lemma}\label{QK} $K\subset Q^*_K$ and
$$
\frac{1}{5}\dist(x_K,S)\le \diam Q_K\le \dist(x_K,S).
$$
\end{lemma}
{\bf Proof.} Since $x_K\in Q_K$, by Theorem \ref{Wcov},
(ii),
$$
\diam Q_K\le \dist(Q_K,S)\le \dist(x_K,S).
$$
Applying again property (ii) of Theorem \ref{Wcov}, we
obtain
$$
\dist(x_K,S)\le \diam Q_K+\dist(Q_K,S)\le 5\diam Q_K.
$$
\par This inequality and \rf{KA} imply
$$
\diam K\le \frac{1}{40}\dist(x_K,S)\le \frac{1}{8}\diam Q.
$$
Since $Q_K\cap K\ne\emp$, we obtain the required inclusion
$K\subset(1+\frac{1}{8})Q_K=:Q^*_K$.\bx
\begin{lemma}\label{KDif}(Brudnyi \cite{Br1})
Let $Q$ be a cube in $\RN$ and let $g\in C^{\infty}(Q)$.
Then for every $1\le u\le\infty$
$$
\Ec_k(g;Q)_{L_u}\le C(\diam Q)^k\max_{|\alpha|=k}
\|D^\alpha g\|_{L_\infty(Q)}.
$$
\end{lemma}
\begin{lemma}\label{Ksmall} For every cube $K$ satisfying
\rf{KA} and every $1\le u\le \infty$ we have
$$
\Ec_k(\tf;K)_{L_u}\le C\left(\frac{\diam K}{\dist
(x_K,S)}\right)^k\max\{\|P_Q-P_{Q_K}\|_{L_\infty(Q)}:~
Q^*\cap K\ne\emp\}.
$$
\end{lemma}
{\bf Proof.} Since $K\subset\RN\setminus S$, $\tf|_K\in
C^\infty(K)$, so that by Lemma \ref{KDif}
$$
\Ec_k(\tf;K)_{L_u}\le C(\diam K)^k\max_{|\alpha|=k}
\|D^\alpha\tf\|_{L_\infty(K)}.
$$
Since $K\subset Q^*_K$, see Lemma \ref{QK}, by properties
of partition of unity and by Leibnitz's formula for every
$|\alpha|=k$ we have
\be \|D^\alpha\tf\|_{L_\infty(K)}&=&
\|D^\alpha(\sum\limits_{Q\in W_S}
\varphi_Q(P_Q-P_{Q_K}))\|_{L_\infty(K)}\nn\\
&\le& C\sum\limits_{Q^*\cap
K\ne\emp}~\sum\limits_{\alpha=\alpha_1+\alpha_2}
\|D^{\alpha_1}\varphi_Q\|_{L_\infty(K)} \cdot
\|D^{\alpha_2}(P_Q-P_{Q_K})\|_{L_\infty(K)}\nn\\
&\le& C\sum\limits_{Q^*\cap
K\ne\emp}~\sum\limits_{\alpha=\alpha_1+\alpha_2} (\diam
Q)^{-|\alpha_1|}
\|D^{\alpha_2}(P_Q-P_{Q_K})\|_{L_\infty(Q^*_K)}.\nn \ee
By Markov's inequality and Proposition \ref{Pol}
\be
\|D^{\alpha_2}(P_Q-P_{Q_K})\|_{L_\infty(Q^*_K)} &\le&
C(\diam
Q^*_K)^{-|\alpha_2|}\|P_Q-P_{Q_K}\|_{L_\infty(Q^*_K)}\nn\\
&\le& C(\diam
Q_K)^{-|\alpha_2|}\|P_Q-P_{Q_K}\|_{L_\infty(Q_K)}\nn
\ee
so that
$$
\Ec_k(\tf;K)_{L_u}\le C\left(\frac{\diam K}{\diam
Q_K}\right)^k\sum\{\|P_Q-P_{Q_K}\|_{L_\infty(Q)}:~ Q^*\cap
K\ne\emp\}.
$$
But by  Lemma \ref{QK} $\diam Q_K\approx \dist(x_K,S)$, and
the result follows.\bx
\par Recall  that $a_{x_K}$ stands for a point nearest
to $x_K$ on $S$. Denote
$$
\tQ_K:= Q(a_{x_K},2\dist(x_K,S)).
$$
Then inequality \rf{KA} immediately implies that
$\tQ_K\supset K$.
\par We put
$$
\Ac_K:=\{Q\in\Qc_2(\tQ_K):~\diam Q\le 4\dist(x_K,S)\}.
$$
\begin{lemma}\label{QW} For every cube $K$ satisfying
\rf{KA} and for every polynomial $P_0\in\Pc_{k-1}$
$$
\Ec_k(\tf;K)^u_{L_u}\le C\left(\frac{\diam
K}{\dist(x_K,S)}\right)^{ku}|\tQ_K|^{-1}\sum_{Q\in\Ac_K}
\|P_Q-P_0\|^u_{L_u(Q)}.
$$
\end{lemma}
{\bf Proof.} For each $Q\in W_S$ such that $Q^*\cap
K\ne\emp$ we have
$$
\|P_Q-P_{Q_K}\|_{L_\infty(Q)}\le \|P_Q-P_0\|_{L_\infty(Q)}
+\|P_{Q_K}-P_0\|_{L_\infty(Q)}.
$$
Since $Q^*_K\supset K$, see Lemma \ref{QK}, we have
$Q^*\cap Q^*_K\ne\emp$ so that by Lemma \ref{Wadd}, (1),
$\diam Q\approx\diam Q_K$. Then by Corollary \ref{CorB}
$$
\|P_{Q_K}-P_0\|_{L_\infty(Q)}\le
C\|P_{Q_K}-P_0\|_{L_\infty(Q_K)}.
$$
In turn, by Corollary \ref{CorA}
$$
\|P_{Q}-P_0\|^u_{L_\infty(Q)}\le
C|Q|^{-1}\|P_{Q}-P_0\|^u_{L_u(Q)} \le
C|\tQ_K|^{-1}\|P_{Q}-P_0\|^u_{L_u(Q)} .
$$
Hence
$$
\max\{\|P_{Q}-P_{Q_K}\|^u_{L_\infty(Q)}:Q^*\cap K\ne\emp\}
\le C|\tQ_K|^{-1} \max\{\|P_{Q}-P_0\|^u_{L_u(Q)}:Q^*\cap
K\ne\emp\}
$$
so that by Lemma \ref{Ksmall}
\be \Ec_k(\tf;K)^u_{L_u}&\le& C\left(\frac{\diam
K}{\dist(x_K,S)}\right)^{ku}
\max\{\|P_Q-P_{Q_K}\|^u_{L_\infty(Q)}:~
Q^*\cap K\ne\emp\}\nn\\
&\le& C\left(\frac{\diam
K}{\dist(x_K,S)}\right)^{ku}|\tQ_K|^{-1}
\max\{\|P_{Q}-P_0\|^u_{L_u(Q)}: ~ Q^*\cap K\ne\emp\}.\nn
\ee
\par Since $K\subset\tQ_K$, by definition of the family
$\Qc_2$, see \rf{AA}, every $Q\in W_S$ such that $Q^*\cap
K\ne\emp$ belongs to $\Qc_2(\tQ_K)$. Moreover, by Lemma
\ref{Wadd}, (1),  $\diam Q\le 4\diam Q_K$ and by Lemma
\ref{QK} $\diam Q_K\le \dist(x_K,S)$. Hence $\diam Q\le
4\dist(x_K,S)$ proving that $Q\in \Ac_K$. This shows that
the latter maximum can be taken over family $\Ac_K$. The
the lemma is proved. \bx
\par We put
\bel{YY}
\QB:= 25\tQ_K=Q(a_{x_K},50\dist(x_K,S)).
\ee
\begin{lemma}\label{ES1} Suppose that a cube $K$ satisfies
\rf{KA} and $\dist(x_K,S)\le\delta_S/4$. Then
$$
\Ec_k(\tf;K)_{L_u}\le C\left(\frac{\diam
K}{\dist(x_K,S)}\right)^{k}\Ec_k(f;\QB)_{L_u(S)}.
$$
\end{lemma}
{\bf Proof.} Since $\dist(x_K,S)\le\delta_S/4$,
\be
\Ac_K&:=&\{Q\in\Qc_2(\tQ_K):~\diam Q\le
4\dist(x_K,S)\}\nn\\
&\subset& \{Q\in\Qc_2(\tQ_K):~\diam Q\le
\delta_S\}=:\Qc_3(\tQ_K),\nn \ee
see \rf{AA3}. Hence by Lemma \ref{QW} for every
$P_0\in\Pc_{k-1}$ we have
$$ \Ec_k(\tf;K)^u_{L_u}\le C\left(\frac{\diam
K}{\dist(x_K,S)}\right)^{ku}|\tQ|^{-1}
\sum_{Q\in\Qc_3(\tQ_K)} \|P_Q-P_0\|^u_{L_u(Q)}. $$
By Lemma \ref{ETWO}
$$
\sum_{Q\in\Qc_3(\tQ_K)} \|P_{Q}-P_0\|^u_{L_u(Q)} \le
C\|f-P_0\|^u_{L_u((25\tQ_K)\cap
S)}=C\|f-P_0\|^u_{L_u(\QB\cap S)}.
$$
\par It remains to put $P_0\in\Pc_{k-1}$ to be a
polynomial of the best approximation of $f$ on $\QB\cap S$
in $L_u$-norm and the lemma follows.\bx
\par The last auxiliary result of the section is the
following
\begin{lemma}\label{NormK} For every cube $K$ satisfying
\rf{KA}
$$
\Ec_k(\tf;K)_{L_u}\le C\left(\frac{\diam
K}{\dist(x_K,S)}\right)^{k}
|\QB|^{-\frac{1}{u}}\|f\|_{L_u(\QB\cap S)}.
$$
\end{lemma}
{\bf Proof.} Recall that $P_Q:=0$ if $\diam Q>\delta_S$ so
that
$$
\sum_{Q\in\Qc_2(\tQ_K)} \|P_{Q}\|^u_{L_u(Q)}=
\sum_{Q\in\Qc_3(\tQ_K)} \|P_{Q}\|^u_{L_u(Q)}.
$$
By Lemma \ref{QW} with $P_0=0$ we obtain
$$
\Ec_k(\tf;K)^u_{L_u}\le C\left(\frac{\diam
K}{\dist(x_K,S)}\right)^{ku}
|\tQ_K|^{-1}\sum_{Q\in\Qc_3(\tQ_K)}
\|P_Q\|^u_{L_u(Q)}.
$$
Hence by Lemma \ref{ETWO} (with $P_0=0$) we have
$$
\Ec_k(\tf;K)^u_{L_u}\le C\left(\frac{\diam
K}{\dist(x_K,S)}\right)^{ku}
|\tQ_K|^{-1}\|f\|^u_{L_u((25\tQ_K)\cap S)}
$$
which implies the lemma because $\QB:=25\tQ_K$.\bx
\par We are in a position to finish the proof of Theorem
\ref{ExtLoc}. Let us fix $x\in\RN$ and $t>0$ and consider
four cases.
\par {\it Case 1.} $80t\le\dist(x,S)$ and
$r^{(x,t)}\le\delta_S$. Recall that
$ r^{(x,t)}:=50\max\{80t,\dist(x,S)\}$
so that in our case $r^{(x,t)}=50\dist(x,S)$. In turn,
$$
K^{(x,t)}:=Q(a_x,r^{(x,t)})=Q(a_x,50\dist(x,S)),
$$
see \rf{RXT} and \rf{KXT}.
\par Put $K:=Q(x,t)$. Then $\diam K=2t$
(recall that we measure distances in the uniform norm) so
that $\diam K\le\dist(x,S)/40$. Moreover,
$$
r^{(x,t)}=50\dist(x,S)\le\delta_S
$$
which, in particular, implies that
$\dist(x,S)\le\delta_S/2$. Thus $K$ satisfies conditions of
Lemma \ref{ES1}. By this lemma
$$
\Ec_k(\tf;K)_{L_u}\le
C\left(\frac{t}{\dist(x,S)}\right)^{k}\Ec_k(f;\QB)_{L_u(S)}
$$
where
$
\QB:=Q(a_{x},50\dist(x,S))=K^{(x,t)},
$
see \rf{YY}. Since $80t\le\dist(x,S)$, we have
$\dist(x,S)^k\approx t^k+\dist(x,S)^k$ proving Theorem
\ref{ExtLoc} in the case under consideration.
\par {\it Case 2.} $80t\le\dist(x,S)$ and
$r^{(x,t)}>\delta_S$.
\par We treat this case in the same way as the previous
one. The only difference is we apply Lemma \ref{NormK}
rather than Lemma \ref{ES1}.
\par {\it Case 3.} $80t>\dist(x,S)$ and
$r^{(x,t)}\le\delta_S$. In this case $r^{(x,t)}=50\cdot
80t=4000t$ so that $4000t\le \delta_S$. Recall that
$\|a_x-x\|_\infty=\dist(x,S)$ so that
$$
K=Q(x,t)\subset Q(a_x,\dist(x,S)+t)\subset Q(a_x,81t).
$$
We put $\OK:=Q(a_x,81t)$ so that $K\subset\OK$. Then by
\rf{ESUB}
\bel{BigK}
\Ec_k(\tf;K)_{L_u}\le C \Ec_k(\tf;\OK)_{L_u}
\ee
and by Theorem \ref{Approx}
$\Ec_k(\tf;\OK)_{L_u}\le C\Ec_k(f;25\OK)_{L_u(S)}. $
Observe that
\bel{KXF} 25\OK\subset
K^{(x,t)}:=Q(a_{x},r^{(x,t)})=Q(a_{x},4000t) \ee
so that by \rf{ESUB}
$
\Ec_k(f;25\OK)_{L_u(S)}\le C \Ec_k(f;K^{(x,t)})_{L_u(S)}.
$
\par It remains to note that $t^k+\dist(x,S)^k\approx t^k$ and
Case 3 is proved.
\par {\it Case 4.} $80t>\dist(x,S)$ and
$r^{(x,t)}>\delta_S$. We preserves notation of the previous
case so that we assume that inequality \rf{BigK} holds.
Clearly,
$
\Ec_k(\tf;\OK)_{L_u}\le |\OK|^{-\frac{1}{u}}
\|\tf\|_{L_u(\OK)}
$
so that by Proposition \ref{ELU}
$$
\Ec_k(\tf;\OK)_{L_u}\le C|\OK|^{-\frac{1}{u}}
\|f\|_{L_u((25\OK)\cap S)}.
$$
Combining this with \rf{BigK} and \rf{KXF} we obtain
$
\Ec_k(\tf;K)_{L_u}\le C|\OK|^{-\frac{1}{u}}
\|f\|_{L_u(K^{(x,t)}\cap S)}.
$
 Since $|\OK|\approx |K^{(x,t)}|$ and
$t^k+\dist(x,S)^k\approx t^k$ this proves Case 4 and the
theorem.\bx
\SECTLONG{4. Estimates of sharp maximal functions: proofs
of}{Theorem \ref{EXT1} and Theorem \ref{EXT2}.}{4}
\indent
\par To formulate the main result of the section we fix
parameters $s\ge 0$, $k\in\N\cup\{0\}$, $0<q\le\infty$,
$1\le u\le\infty,$ and put $\vv:=(s,k,q,u)$.
Given a function $f\in L_{u,\,loc}(S)$, we let $\fs$ denote
a {\it generalized sharp maximal function} of $f$ on $S$:
\bel{MFG} \fs(x):=\left\{\int_0^\infty\left(
\frac{\Ec_k(f;Q(x,t))_{L_u(S)}}{t^s}\right)^q\frac{dt}{t}
\right\}^\frac{1}{q},\ \ \ \ \ x\in S, \ee
if $q<\infty$, and
$$
\fs(x):=\sup_{t>0} \frac{\Ec_k(f;Q(x,t))_{L_u(S)}}{t^s},\ \
\ \ \ x\in S,
$$
if $q=\infty$. We write $f^\sharp_{\vv}$ for
$f^\sharp_{\vv,\RN}$.
\par As usual we put
$M_uf(x):=(M(|f|^u)(x))^{\frac{1}{u}}$ where $M$ is the
Hardy-Littlewood maximal function
$$
Mf(x):=\sup_{t>0}\frac{1}{|Q(x,t)|}\int_{Q(x,t)}|f(y)|dy.
$$
We recall that by the Hardy-Littlewood-Wiener maximal
inequality, see e.g. \cite{St}, for every $0<u<p\le\infty$
and $g\in L_p(\RN)$
\bel{HLW}
\|M_ug\|_{L_p(\RN)}\le C\|g\|_{L_p(\RN)}.
\ee
\begin{theorem}\label{MF} Let $S$ be a regular subset of
$\RN$ and let $f\in L_{u,\,loc}(S)$. Assume that $1\le
u\le\infty$,  and  $0\le s<k$ if $0<q\le\infty$ or $0\le
s\le k$ if $q=\infty$. Then
$$
\tfs(x)\le C\{M((\fs)^\cw)(x)+M_u(f^\cw)(x)\}, \ \ \ \ \
x\in\RN.
$$
\end{theorem}
\par Recall that $\tf$ stands for the extension of $f$
defined by formula \rf{ExtOp}. Recall also that $h^\cw$
where $h$ is a function on $S$ denotes the extension by $0$
of $h$ from $S$ on all of $\RN$.
{\bf Proof.} We will prove the result for $0<q<\infty$; the
reader can easy modify this proof for the case $q=\infty$.
Let us consider two cases.
\par {\it Case 1.} We assume that
\bel{Case1} \1\le \delta_S/50. \ee
\par Put $\Delta:=\delta_S/4000$. Then
$$
\tfs(x):=\left\{\int_0^\infty\left(
\frac{\Ec_k(\tf;Q(x,t))_{L_u}}{t^s}\right)^q\frac{dt}{t}
\right\}^\frac{1}{q}\le C(I_1+I_2)
$$
where
$$
I_1:=\left\{\int_0^\Delta\left(
\frac{\Ec_k(\tf;Q(x,t))_{L_u}}{t^s}\right)^q\frac{dt}{t}
\right\}^\frac{1}{q} ~ {\rm and} ~
I_2:=\left\{\int_\Delta^\infty\left(
\frac{\Ec_k(\tf;Q(x,t))_{L_u}}{t^s}\right)^q\frac{dt}{t}
\right\}^\frac{1}{q}.
$$
\par Let us estimate $I_1$. We observe that for every
$0<t\le\Delta$ by inequality \rf{Case1} the quantity
$r^{(x,t)}:=50\max\{80t,\1\}$ satisfies the inequality
$r^{(x,t)}\le\delta_S$. Therefore by Theorem \ref{ExtLoc}
$$
\Ec_k(\tf;Q(x,t))_{L_u}\le C\frac{t^k}{t^k+\dist(x,S)^k}
\,\Ec_k(f;K^{(x,t)})_{L_u(S)}.
$$
Recall that $K^{(x,t)}:=Q(a_x,r^{(x,t)})$ where $a_x$ is a
point on $S$ such that
\bel{AX}
\|x-a_x\|=\1.
\ee
Hence
$$
I_1\le\left\{\int_0^\Delta
\left(\frac{t^{k-s}}{t^k+\dist(x,S)^k}\,
\Ec_k(f;K^{(x,t)})_{L_u(S)}\right)^q\frac{dt}{t}
\right\}^\frac{1}{q} \le C(J_1+J_2)
$$
where
$$
J_1:=\left\{\int\limits_0^{\1/80}
\left(\frac{t^{k-s}}{\dist(x,S)^k}\,
\Ec_k(f;K^{(x,t)})_{L_u(S)}\right)^q\frac{dt}{t}
\right\}^\frac{1}{q}
$$
and
$$
J_2:=\left\{\int\limits_{\1/80}^\infty
\left(\frac{\Ec_k(f;K^{(x,t)})_{L_u(S)}}{t^s}\right)^q
\frac{dt}{t}\right\}^\frac{1}{q}.
$$
\par Prove that $J_1\le C\,J_2$. We observe that for every
$0<t\le\1/80$ we have $r^{(x,t)}:=50\1$ so that
$K^{(x,t)}:=Q(a_x,50\1)$. Hence
$$
J_1\le
C\,\frac{\Ec_k(f;Q(a_x,50\1))_{L_u(S)}}{\dist(x,S)^k}
\left\{\int\limits_0^{\1/80} t^{(k-s)q}\frac{dt}{t}
\right\}^\frac{1}{q}.
$$
Since $k>s\ge 0$ or $k\ge s\ge 0$ if $q=\infty$, the latter
integral is equivalent $\1^{(k-s)q}$. Hence
$$ J_1\le
C\,\frac{\Ec_k(f;Q(a_x,50\1))_{L_u(S)}}{\dist(x,S)^s}. $$
By \rf{ESUB} for every $t$ such that $\1<t\le2\1$ we have
$$
\Ec_k(f;Q(a_x,50\1))_{L_u(S)}\approx
\Ec_k(f;Q(a_x,50t))_{L_u(S)}
$$
so that
$$
\frac{\Ec_k(f;Q(a_x,50\1))_{L_u(S)}}{\dist(x,S)^s} \approx
\left\{\int\limits_{\1}^{2\1}
\left(\frac{\Ec_k(f;Q(a_x,50t))_{L_u(S)}}{t^s}\right)^q
\frac{dt}{t} \right\}^\frac{1}{q}.
$$
Observe that for $t\ge\1/80$ we have
$r^{(x,t)}:=50\max\{80t,\1\}=4000t$ and
$K^{(x,t)}:=Q(a_x,4000t)$ so that
\bel{J2E}
J_2\approx \left\{\int\limits_{\1}^\infty
\left(\frac{\Ec_k(f;Q(a_x,50t))_{L_u(S)}}{t^s}\right)^q
\frac{dt}{t}\right\}^\frac{1}{q}
\ee
proving the required inequality $J_1\le C J_2$.
\par Let us estimate $J_2$. To this end we put
$\tK:=Q(x,2\1)$. Prove that for each $y\in\tK\cap S$ we
have $J_2\le C\fs(y).$ In fact, by \rf{AX}
$$
\|y-a_x\|\le \|y-x\|+\|x-a_x\|\le 3\1.
$$
Hence for every $t>50\1$ we have $Q(a_x,t)\subset Q(y,2t)$.
This inclusion and \rf{ESUB} imply
$$
\Ec_k(f;Q(a_x,t))_{L_u(S)}\le C\Ec_k(f;Q(y,2t))_{L_u(S)}.
$$
so that by \rf{J2E}
$$
J_2\le C \left\{\int\limits_{\1}^\infty
\left(\frac{\Ec_k(f;Q(y,100t))_{L_u(S)}}{t^s}\right)^q
\frac{dt}{t}\right\}^\frac{1}{q}\le C \fs(y)
$$
proving the required inequality $J_2\le C\fs(y).$ By this
inequality
$$
J_2\le C\frac{1}{|\tK\cap S|}\int\limits_{\tK\cap S}
\fs(y)dy.
$$
Since $\1\ge\delta_S/50$, see \rf{Case1}, by \rf{AX} we
have
$$
Q(a_x,\1)\subset Q(x,2\1)=:\tK.
$$
Since $S$ is regular and $\1\le\delta_S$,
$$
|\tK\cap S|\ge |Q(a_x,\1)\cap S|\ge
\theta_S|Q(a_x,\1)|\approx |\tK|
$$
so that $ |\tK\cap S|\approx |\tK| $. Hence
$$
J_2\le C\frac{1}{|\tK|}\int\limits_{\tK\cap S} \fs(y)dy \le
C M(\fs)^\cw(x).
$$
Combining this with the estimate $J_1\le C J_2$ we conclude
that $I_1\le C M(\fs)^\cw(x)$.
\par Let us prove that $I_2\le C M_u(f^\cw)(x)$.
We recall that $\1\le\delta_S/50$ so that for every
$t>\Delta :=\delta_S/4000$
$$
r^{(x,t)}:=50\max\{80t,\1\}=4000t>\delta_S
$$
and $K^{(x,t)}:=Q(a_x,r^{(x,t)})=Q(a_x,4000t).$
Therefore by Theorem \ref{ExtLoc} and \rf{K0}
$$
\Ec_k(\tf;Q(x,t))_{L_u}\le
C\left(\frac{1}{|K^{(x,t)}|}\int_{K^{(x,t)}\cap S}
|f|^udy\right)^{\frac{1}{u}}.
$$
We put $\OK:=Q(x,4080t)$. Since $\1\le\delta_S/50\le 80t$,
by \rf{AX} we have
$$
K^{(x,t)}\subset Q(x,\1+4000t)\subset\OK.
$$
Moreover, we also obtain an equivalence
$|K^{(x,t)}|\approx|\OK|$. Hence
$$
\Ec_k(\tf;Q(x,t))_{L_u}\le
C\left(\frac{1}{|\OK|}\int_{\OK\cap S}
|f|^udy\right)^{\frac{1}{u}}\le CM_u(f^\cw)(x).
$$
This implies
$$
I_2:=\left\{\int_\Delta^\infty\left(
\frac{\Ec_k(\tf;Q(x,t))_{L_u}}{t^s}\right)^q\frac{dt}{t}
\right\}^\frac{1}{q}\le
C(M_u(f^\cw)(x))\left(\int_\Delta^\infty
t^{-sq-1}dt\right)^{\frac{1}{q}}
$$
proving the required $I_2\le C M_u(f^\cw)(x)$. Finally, we
obtain
$$
\tfs(x)\le C\{I_1+I_2\}\le
C\{M((\fs)^\cw)(x)+M_u(f^\cw)(x)\}
$$
which proves the theorem in the first case.
\par {\it Case 2.} $\1>\delta_S/50$. In this case
$r^{(x,t)}:=50\max\{80t,\1\}>\delta_S$ so that by Theorem
\ref{ExtLoc} and \rf{K0}
\bel{EK1}
\Ec_k(\tf;Q(x,t))_{L_u}\le
C\frac{t^k}{t^k+\dist(x,S)^k}\,
|K^{(x,t)}|^{-\frac{1}{u}}\|f\|_{L_u(K^{(x,t)}\cap S)}.
\ee
Recall that $K^{(x,t)}:=Q(a_x,r^{(x,t)})$.
\par Put $K':=Q(x,2r^{(x,t)})$. Clearly,
$r^{(x,t)}\ge 50\1\ge\1$ so that by \rf{AX}
$$
K^{(x,t)}:=Q(a_x,r^{(x,t)})\subset Q(x,\1+r^{(x,t)})\subset
Q(x,2r^{(x,t)})=:K'.
$$
Hence
\bel{EK2}
|K^{(x,t)}|^{-\frac{1}{u}}\|f\|_{L_u(K^{(x,t)}\cap S)}\le
C|K'|^{-\frac{1}{u}}\|f\|_{L_u(K'\cap S)}\le
CM_u(f^\cw)(x).
\ee
Estimates \rf{EK1} and \rf{EK2} and definition \rf{MFG}
imply
$$
\tfs(x)\le CM_u(f^\cw)(x)
\left(\int\limits_0^\infty\frac{t^{(k-s)q}}
{(t^k+\dist(x,S)^k)^q}\frac{dt}{t}\right)^\frac{1}{q}. $$
Since $\1\ge\delta_S/50$ and $k>s$ (or $k\ge s$ if
$q=\infty$), the latter integral is bounded by a constant
depending only on $s,k,q$ and $\delta_S$. This proves that
in the case under consideration $\tfs(x)\le
CM_u(f^\cw)(x)$.
\par Theorem \ref{MF} is completely proved.\bx
\par Let us formulate a corollary of this result. To this
end we introduce a slight ge\-ne\-ralization of the maximal
function \rf{MFG}: given $\vv=(s,k,q,u)$ and
$0<\Delta\le\infty$, we put
$$
\fsd(x):=\left\{\int\limits_0^{\Delta}
\left(\frac{\Ec_k(f;Q(x;t))_{L_u(S)}}{t^s}
\right)^q\frac{dt}{t}\right\}^\frac{1}{q},\ \ \ \ \ x\in S,
$$
(with the standard modification for $q=\infty$).
\begin{theorem}\label{MLP}  Suppose that $1\le
u<p\le\infty$, $0<q\le\infty$ and  $k>s\ge 0$ or $k\ge s\ge
0$ if $q=\infty$. If $S$ is a regular set and $f\in
L_p(S)$, then
$$
\|\tfsd\|_{L_p(\RN)}\le
C(\|\fsd\|_{L_p(S)}+\|f\|_{L_p(S)}).
$$
Here the constant $C$ depends also on $\Delta$.
\end{theorem}
{\bf Proof.} Clearly, $\tfsd\le
\tfs\,(:=(\tf)^\sharp_{\vv,\infty,\RN})$ so that
\bel{D1}
\|\tfsd\|_{L_p(\RN)}\le\|\tfs\|_{L_p(\RN)}.
\ee
By Theorem \ref{MF}
$$
\|\tfs\|_{L_p(\RN)}\le C(\|M(\fs)^\cw\|_{L_p(\RN)}
+\|M_u(f^\cw)\|_{L_p(\RN)})
$$
so that by \rf{HLW}
\bel{D2}
\|\tfs\|_{L_p(\RN)}\le C(\|\fs\|_{L_p(S)}
+\|f\|_{L_p(S)}).
\ee
On the other hand, for every $x\in\RN$,
$$
\fs(x):= \left\{\int_0^\infty\left(
\frac{\Ec_k(f;Q(x,t))_{L_u(S)}}{t^s}\right)^q\frac{dt}{t}
\right\}^\frac{1}{q}\le C(\fsd(x)+J(x))
$$
where
$$
J(x):= \left\{\int_\Delta^\infty\left(
\frac{\Ec_k(f;Q(x,t))_{L_u(S)}}{t^s}\right)^q\frac{dt}{t}
\right\}^\frac{1}{q}.
$$
Observe that
$$
\Ec_k(f;Q(x,t))_{L_u(S)}\le\left(\frac{1}{|Q(x,t)|}
\int_{Q(x,t)\cap S}|f|^udy\right)^{\frac{1}{u}}\le
M_u(f^\cw)(x).
$$
Hence
$$
J(x)\le CM_u(f^\cw)(x)\left\{\int_\Delta^\infty t^{-sq-1}dt
\right\}^\frac{1}{q}\le C \Delta^{-sq}M_u(f^\cw)(x).
$$
Thus $\fs(x)\le C(\fsd(x)+M_u(f^\cw)(x))$ so that
$$
\|\fs\|_{L_p(S)}\le C(\|\fsd\|_{L_p(S)}
+\|M_u(f^\cw)\|_{L_p(S)})\le C(\|\fsd\|_{L_p(S)}
+\|f\|_{L_p(S)}).
$$
This inequality, \rf{D1} and \rf{D2} imply the statement of
the theorem.\bx
\par {P\,r\,o\,o\,f\,s} of Theorems \ref{EXT1}
and \ref{EXT2}.
\par Observe that for every locally integrable extension $F$
of $f$ on all of $\RN$ and for each cube $Q$ centered in
$S$ we have $\Ec_k(f;Q)_{L_1(S)}\le \Ec_k(F;Q)_{L_1}$ so
that $f^\sharp_{k,S}\le F^\sharp_{k}$ on $S$.
Then by \rf{AN}
$$
\|f\|_{L_p(S)}+\|f^\sharp_{k,S}\|_{L_p(S)}\le
\|F\|_{L_p(\RN)}+\|F^\sharp_{k}\|_{L_p(\RN)} \le
C\|F\|_{\WKP}
$$
proving that
$$
\|f\|_{L_p(S)}+\|f^\sharp_{k,S}\|_{L_p(S)}\le
C\|f\|_{\WKP|_S}.
$$
In a similar way using equivalence \rf{Crit} we show that
$$
\|f\|_{L^p(S)}+ \left\|\left(\int_0^1
\left(\frac{\Ec_k(f;Q(\cdot,t))_{L_1(S)}}{t^s}\right)^q\,
\frac{dt}{t}\right)^{\frac{1}{q}}\right\|_{L_p(S)}\le
C\|f\|_{\FS|_S}.
$$
\par To prove the opposite inequalities
we observe that by Proposition \ref{ELU} we have
$
\|\tf\|_{L_p(\RN)}\le C\|f\|_{L_p(S)}.
$
Moreover, by Theorem \ref{MLP} with $1<p\le\infty, k=s,
q=\infty,u=1$ and $\Delta=\infty$
$$
\|(\tf)^\sharp_k\|_{L^p(\RN)}\le
C(\|f^\sharp_{k,S}\|_{L^p(S)}+\|f\|_{L_p(S)}).
$$
Hence
$$
\|f\|_{\WKP|_S}\le \|\tf\|_{\WKP}\le
C(\|\tf\|_{L_p(\RN)}+\|(\tf)_{k}^{\sharp}\|_{L_p(\RN)}) \le
C(\|f\|_{L_p(S)}+\|f^\sharp_{k,S}\|_{L^p(S)})
$$
proving \rf{Sob}. In a similar way we prove equivalence
\rf{FSnorm} applying Theorem \ref{MLP} with $0<s<k,$
$1<p\le\infty,$ $1\le q\le\infty,u=1$ and $\Delta=1$.\bx
\SECT{5.  Besov spaces on regular subsets of $\RN$.}{5}
\par We turn to the problem of an intrinsic
characterization of traces of the Besov spaces to regular
subsets of $\RN$. First we recall one of the equivalent
definitions of the Besov spaces: a function $f$ defined on
$\RN$ belongs to the space $\BS,$ $1\le p\le\infty,
0<q\le\infty,$ $s>0,$ if $f\in L_p(\RN)$ and its modulus of
continuity of order $k$ in $L_p$
$$
\omega_k(f;t)_{L_p}:=\sup_{\|h\|\le t}\|\triangle^k_h
f\|_{L_p(\RN)}
$$
satisfies the inequality
$$\int_0^1 \left(\frac{\omega_k(f;t)_{L_p}}{t^s}\right)^q\,
\frac{dt}{t}<\infty $$
($\sup\limits_{0<t\le 1}t^{-s}\omega_k(f,t)_{L_p}<\infty$
if $q=\infty$). Here $k>s$ is an arbitrary integer and as
usual given $x,h\in\RN$,
$$
\triangle^k_h f(x):=\sum_{j=0}^k(-1)^{k-j}{n\choose
j}f(x+jh).
$$
$\BS$ is normed by
\bel{BN}
\|f\|_{\BS}:=\|f\|_{L_p(\RN)}+\left(\int_0^1
\left(\frac{\omega_k(f;t)_{L_p}}{t^s}\right)^q\,
\frac{dt}{t}\right)^{\frac{1}{q}}
\ee
(modification if $q=\infty$).
\par Similar to the case of Sobolev and $F$-spaces the main
point of our approach to intrinsic characterization of
traces of Besov spaces is local approximations theory.
\par As we have mentioned above this theory gives
a unified approach to various types of function spaces
based on the concept of local best approximation by
polynomials, see definitions \rf{ESP} and \rf{EQ}.
Comparing classical approximation theory and local
approximation theory we observe that one basic goal of
classical approximation theory is to study functions via
the behavior of their best approximations  as a function of
the {\it degree} of the approximating polynomials on a {\it
fixed} set. In local approximation theory we have a similar
goal, but rather than doing all approximations on a fixed
set, we do it on a {\it variable} cube. We can think of it
as a ''window'' which we can slide around, enlarge and
contract, ''looking'' through it at the function's graph.
Each time we consider approximation on the cube by
polynomials of a {\it fixed} (maybe small) degree, and we
study the behavior of the best approximations as a function
of the position and size of the sliding cube.
\par As an important example illustrating this idea we
present so-called an ``atomic" decomposition of the modulus
of continuity due to Brudnyi \cite{Br1,Br2}, see also
\cite{Br3,Br6,Br4}. This basic fact of local approximation
theory states that for every $0<p\le\infty$, $k\in\N$, and
every function $f\in L_{p,\,loc}(\RN)$
\bel{AD}
\omega_k(f;t)_{L_p}\approx\sup_{\pi}\left\{\sum_{Q\in \pi}
E_k(f;Q)^p_{L_p}\right\}^{\frac{1}{p}}
\ee
where the supremum is taken over all {\it packings} $\pi$
of equal cubes in $\RN$ with diameter $t$. (Hereafter
``packing" means a {\it finite family of disjoint cubes} in
$\RN$.) Observe that equivalence \rf{AD} remains true if
$\pi$ runs over all packings of equal cubes with diameter
{\it at most} $t$, see \cite{Br1}.
\par This result motivates the following definition,
see \cite{Br4}: given $k\in\N$, $0<u,p\le\infty$, and a
function $f\in L_{u,\,loc}$, by
$\Omega_{k,p}(f;\cdot)_{L_u}$ we denote the {\it
$(k,p)$-modulus of continuity of $f$ in $L_u$,} i.e., a
function of $t>0$ defined by the following formula
\bel{KPU}
\Omega_{k,p}(f;t)_{L_u}:=\sup_{\pi}\left\{\sum_{Q\in \pi}
|Q|\,\Ec_k(f;Q)^p_{L_u}\right\}^{\frac{1}{p}}. \ee
Here $\pi$ runs over all packings of equal cubes in $\RN$
with diameter $t$. (This definition is a slight
modification of that given in \cite{Br1} where the supremum
is taken over all packings $\pi$ of equal cubes with
diameter {\it at most} $t$.)
\par We note two important properties of the
$(k,p)$-modulus of continuity. First of them  is the
following equivalence, see \cite{Br4}, Chapter 3, and Lemma
\ref{INR}:
\bel{IN1}
\Omega_{k,p}(f;t)_{L_u}\approx
\|\Ec_k(f;Q(\cdot;t))_{L_u}\|_{L_p(\RN)}.
\ee
\par In particular,  from \rf{IN1} and \rf{AD} it follows
that
\bel{EQM} \omega_k(f;t)_{L_p}\approx
\Omega_{k,p}(f;t)_{L_p}\approx
\|\Ec_k(f;Q(\cdot;t))_{L_p}\|_{L_p(\RN)}, \ \ \ \ \ t>0.
\ee
\par The second property clarifies connections between
the $(k,p)$-moduli of continuity in different metrics.
Clearly,
$\Omega_{k,p}(f;\cdot)_{L_u}\le\Omega_{k,p}(f;t)_{L_p}$
whenever $0<u\le p$. On the other hand, Brudnyi
\cite{Br4,Br1} has proved that for every $1\le u\le p$
\bel{DIFM} \Omega_{k,p}(f;t)_{L_p}\le
C\int_0^t\frac{\Omega_{k,p}(f;\tau)_{L_u}}{\tau}\,d\tau,
 \ \ \ \ \ t>0. \ee
\par Now combining definition \rf{BN}, equivalence \rf{IN1}
and inequality \rf{DIFM} and applying the Hardy inequality
we obtain characterization \rf{ChB} of Besov functions on
$\RN$ via local approximations.
\par Let us generalize  definition \rf{KPU} for the case of
a measurable subset $S\subset\RN$ and a function $f\in
L_{u,\,loc}(S)$. We define the $(k,p)$-modulus of
continuity of $f$ in $L_u(S)$ (\cite{Br4}) by letting
\bel{OKP}
\Omega_{k,p}(f;t)_{L_u(S)}:=\sup_{\pi}\left\{\sum_{Q\in
\pi} |Q\cap
S|\,\Ec_k(f;Q)^p_{L_u(S)}\right\}^{\frac{1}{p}}. \ee
Here $\pi$ runs over all packings of equal cubes centered
in $S$ with diameter $t$.
\par Let us show that an analog of equivalence \rf{IN1} is
true for $ \Omega_{k,p}(f;\cdot)_{L_u(S)}$ as well. To
prove this we need the following simple combinatorial
lemma.
\begin{lemma}\label{FQ} Let $\pi$ be a family of equal
cubes such that
$ \sum\{\chi_Q:~Q\in\pi\}\le l $
where $l$ is a positive constant. Then a family of cubes
$\{2 Q:~Q\in\pi\}$ can be represented as union of at most
$m=m(l)$ packings.
\end{lemma}
\par In particular, from the lemma and definition \rf{OKP}
it easily follows that $\Omega_{k,p}(f;\cdot)_{L_u(S)}$ is
a quasi-monotone function, i.e.,
\bel{QMON} \Omega_{k,p}(f;t)_{L_u(S)}\le C\,
\Omega_{k,p}(f;2t)_{L_u(S)}, \ \ \ \ \ t>0. \ee
\begin{lemma}\label{INR} Let $0<u,p\le\infty$ and $k\in\N$.
Then for every function $f\in L_{u,\,loc}(S)$
$$
\frac{1}{C}\,\Omega_{k,p}(f;t/4)_{L_u(S)}\le
\|\Ec_k(f;Q(\cdot,t))_{L_u(S)}\|_{L_p(S)} \le
C\,\Omega_{k,p}(f;t)_{L_u(S)},\ \ \ \ \ t>0.
$$
\end{lemma}
{\bf Proof.} We will mainly follow a scheme of the proof
given in \cite{Br4} for the case $S=\RN$. Fix $t>0$ and
consider a packing $\pi$ of equal cubes with diameter $t$
centered in $S$. Then for each $Q\in\pi$ and every $x\in
Q\cap S$ we have $Q\subset Q(x,4t)$ so that by \rf{ESUB}
$$
\Ec_k(f;Q)_{L_u(S)} \le C\Ec_k(f;Q(x,4t))_{L_u(S)}.
$$
Hence
$$
|Q\cap S|\,\Ec_k(f;Q)^p_{L_u(S)}\le C\int_{Q\cap S}
\Ec_k(f;Q(x,4t))^p_{L_u(S)}dx,\ \ \ \ \ x\in Q\cap S.
$$
Thus
\be \sum_{Q\in\pi} |Q\cap S|\,\Ec_k(f;Q)^p_{L_u(S)} &\le&
 C\sum_{Q\in\pi}\int_{Q\cap S}
 \Ec_k(f;Q(x,4t))^p_{L_u(S)}dx\nn\\
&\le& C\int_{S}\Ec_k(f;Q(x,4t))^p_{L_u(S)}dx\nn
\ee
proving the first inequality of the lemma.
\par To prove the second inequality given $t>0$, we let
$\tilde{\pi}$ denote a covering of $S$ by equal cubes
centered in $S$ with diameter $t/2$ such that
$\sum\{\chi_Q:~Q\in\tilde{\pi}\}\le C(n).$ (The existence
of $\tilde{\pi}$ immediately follows, for instance, from
the Besicovitch theorem, see e.g. Gusman \cite{G}.) Then
$$
\int_{S}\Ec_k(f;Q(x,t))^p_{L_u(S)}dx\le
\sum_{Q\in\tilde{\pi}}\int_{Q\cap
S}\Ec_k(f;Q(x,t))^p_{L_u(S)}dx.
$$
Clearly, for every $Q\in\tilde{\pi}$ and every $x\in Q\cap
S$ we have $Q(x,t)\subset 2Q=Q(x_Q,2t)$ so that by
\rf{ESUB}
$
\Ec_k(f;Q(x,t))_{L_u(S)} \le C\Ec_k(f;2Q)_{L_u(S)}.
$
Hence
$$
\int_{Q\cap S}\Ec_k(f;Q(x,t))^p_{L_u(S)}dx
 \le C|2Q\cap S|\,\Ec_k(f;2Q)^p_{L_u(S)}
$$
so that
$$
\int_{S}\Ec_k(f;Q(x,t))^p_{L_u(S)}dx\le
 C\sum_{Q\in\tilde{\pi}}|2Q\cap S|\,\Ec_k(f;2Q)^p_{L_u(S)}.
$$
By Lemma \ref{FQ} a family of cubes
$\pi:=\{2Q:~Q\in\tilde{\pi}\}$ can be represented in the
form $\pi=\cup\{\pi_i:~i=1,...,m\}$ where $m=m(n)$ and
every family $\pi_i$ is a packing. Hence
$$
\int_{S}\Ec_k(f;Q(x,t))^p_{L_u(S)}dx\le
 C\sum_{i=1}^m\sum_{Q\in\pi_i}
 |Q\cap S|\,\Ec_k(f;Q)^p_{L_u(S)}
 \le Cm\,\Omega_{k,p}(f;t)_{L_u(S)}
$$
proving the lemma.\bx
\par The main result of the section is the following
\begin{theorem}\label{KPS}
Let $S$ be a regular set and let $1\le u\le p\le\infty$.
Then for every function $f\in L_{u,loc}(S)$ and every
$0<t\le 1$
$$
\|\Ec_k(\tf;Q(\cdot,t))_{L_u}\|_{L_p(\RN)}\le C\,
t^k\left\{ \left( \int_t^{1}\left
(\frac{\|\Ec_k(f;Q(\cdot,\tau))_{L_u(S)}\|_{L_p(S)}}
{\tau^k}\right)^p \frac{d\tau}{\tau}\right)^\frac{1}{p}+
\|f\|_{L_p(S)}\right\}
$$
Here $\tf=\Ekus$ is the extension operator defined by
formula \rf{ExtOp}.
\end{theorem}
{\bf Proof.} By Lemma \ref{INR} it is sufficient to show
that for every $0<t\le 1$
\bel{OM} \Omega_{k,p}(\tf;t)_{L_u}\le C\, t^k\left\{ \left(
\int_{t/4}^{1/4}\left (\frac{\Omega_{k,p}(f;\tau)_{L_u(S)}}
{\tau^k}\right)^p \frac{d\tau}{\tau}\right)^\frac{1}{p}+
\|f\|_{L_p(S)}\right\}. \ee
\par First let us estimate $\Omega_{k,p}(f;t/2)_{L_u}$ for
$0<t\le\delta_S/4000$. Fix a family $\pi$ of equal cubes in
$\RN$ of diameter $t/2$. (Thus $Q=Q(x_Q,t)$ for every
$Q\in\pi$.) We let $m=m(t)$ denote a positive integer such
that
\bel{TD}
\frac{1}{2^{m+1}}\,\frac{\delta_S}{50}<80t\le
\frac{1}{2^m}\,\frac{\delta_S}{50}.
\ee
Then for each integer $i,~i<m$ we put
\bel{DPI} \pi_i:=\left\{Q\in\pi:~
\frac{1}{2^{i+1}}\,\frac{\delta_S}{50}<\dist(x_Q,S) \le
\frac{1}{2^i}\,\frac{\delta_S}{50}\right\}. \ee
We also set
\bel{PM} \pi_m:=\left\{Q\in\pi:~ \dist(x_Q,S)
\le\frac{1}{2^m}\,\frac{\delta_S}{50}\right\}. \ee
\par Now following formulas \rf{RXT} and \rf{KXT}
we assign every $Q=Q(x_Q,t)\in\pi$ a number
$$
r^{(x_Q,t)}:=50\max(80t,\dist(x_Q,S))
$$
and a cube $K^{(x_Q,t)}:=Q(a_{x_Q},r^{(x_Q,t)}).$
(Recall that $a_{x_Q}\in S$ and satisfies the equality
$\|x_Q-a_{x_Q}\|_\infty=\dist(x_Q,S)$.)
\par In particular, for every $Q\in\pi_i,\, 0\le i< m,$
we have $\dist(x_Q,S)>80t$, so that in this case
\bel{IR1} r^{(x_Q,t)}=50\dist(x_Q,S) \ee
and $K^{(x_Q,t)}=Q(a_{x_Q},50\dist(x_Q,S))$.
In turn, for $Q\in\pi_m$ we have
$$
\dist(x_Q,S)\le 160t\ \ \ {\rm and}\ \ \ \
r^{(x_Q,t)}\approx t.
$$
Moreover, for every $Q=Q(x_Q,t)\in\pi_i,\, 0\le i\le m,$
\bel{RIM} r^{(x_Q,t)}\le\delta_S. \ee
Observe also that for each $Q=Q(x_Q,t)\in\pi_i$ with $i<0$
\bel{IM} r^{(x_Q,t)}=50\dist(x_Q,S)>\delta_S. \ee
\par We put
\bel{OI}
\Omega_i:=\sum_{Q\in \pi_i}
|Q|\,\Ec_k(\tf;Q)^p_{L_u}.
\ee
\par Let us estimate $\Omega_i$ for $0\le i\le m$. By Theorem
\ref{ExtLoc}, \rf{IR1} and \rf{RIM}, for every $Q\in\pi_i,$
$0\le i< m,$ we have
\bel{ST1} \Ec_k(\tf;Q)_{L_u}\le C\frac{t^k}{\dist(x,S)^k}
\,\Ec_k(f;K^{(x_Q,t)})_{L_u(S)} \ee
where by $K^{(x_Q,t)}:=Q(a_{x_Q},50\dist(x_Q,S)).$
We put
$$
r_i:=2^{-i}\delta_S\ \ \ \ \ {\rm and}\ \ \ \ \
K^{\{Q\}}:=Q(a_{x_Q},r_i).
$$
Since $Q\in\pi_i$, by \rf{DPI} $K^{(x_Q,t)}\subset
K^{\{Q\}}$ and $r_i\approx \dist(x_Q,S)$ so that by
\rf{ESUB} and by \rf{ST1}
\bel{EK}
\Ec_k(\tf;Q)_{L_u}\le C\frac{t^k}{r_i^k}
\,\Ec_k(f;K^{\{Q\}})_{L_u(S)}.
\ee
It can be also readily seen that Theorem \ref{ExtLoc} and
\rf{PM} imply the same estimate for $i=m$ as well. Thus in
what follows we will assume that inequality \rf{EK} is true
for all $i=0,...,m$. Observe also that for each
$Q\in\pi_i$, $i=0,...,m$, we have $Q\subset K^{\{Q\}}.$
\par Now fix $i,\, 0\le i\le m,$ and put
$ \tilde{\pi}_i:=\{K^{\{Q\}}:~Q\in\pi_i\}. $
\par Then by Besicovitch's theorem there is a subfamily
$\pi'_i\subset\tilde{\pi}_i$ such that:
\par (a) for every $K\in\tilde{\pi}_i$ there is a cube
$K'\in\pi'_i$ such that $x_K\in K'$;
\par (b) $\sum\{\chi_{K'}:~K'\in\pi'_i\}\le l(n)$.
\par Now for every cube $K'\in\pi'_i$ we put
$$ \Ac_{K'}:=\{K\in\tilde{\pi}_i:~x_K\in K'\}. $$
Since $\diam K=\diam K'$ and $x_K\in K'$ for every
$K\in\Ac_{K'}$, we have $K\subset 2K'$. Recall also that
$Q\subset K^{\{Q\}}$ so that
\bel{UK}
\cup\{Q:~K^{\{Q\}}\in\Ac_{K'}\}\subset 2K'.
\ee
\par By property (b) of $\pi'_i$ and by Lemma \ref{FQ}
later on we may assume that the family of cubes
$\{2K':~K'\in\pi'_i\}$ is a {\it packing}. By \rf{ESUB} for
every $K\in\Ac_{K'}$ we have
$$\Ec_k(f;K)_{L_u(S)}\le C\Ec_k(f;2K')_{L_u(S)}. $$
Combining this with \rf{EK} we obtain the following
estimate of $\Omega_i$, see \rf{OI}:
$$
\Omega_i\le C\frac{t^{kp}}{r_i^{kp}}
\sum_{Q\in\pi_i}|Q|\Ec_k(f;K^{\{Q\}})^p_{L_u(S)}\le
C\frac{t^{kp}}{r_i^{kp}} \sum_{K'\in\pi'_i}
\left(\sum_{K^{\{Q\}}\in\Ac_{K'}}|Q|\right)
\Ec_k(f;2K')^p_{L_u(S)}.
$$
Hence by \rf{UK}
\bel{ST2} \Omega_i\le C\frac{t^{kp}}{r_i^{kp}}
\sum_{K'\in\pi'_i} |2K'|\, \Ec_k(f;2K')^p_{L_u(S)}. \ee
Recall that
$\diam K'=2r_i:=2^{-i+1}\delta_S\le 2\delta_S $ for every
cube $K'\in\pi'_i$ which implies
$\diam(\frac{1}{2}K')=\frac{1}{2}\diam K'\le\delta_S$.
Since $S$ is regular, we obtain
$$
|2K'|=4^n|(1/2)K'|\le 4^n\delta_S|(1/2)K')\cap S|
 \le 4^n\delta_S|(2K')\cap S|.
$$
Hence
$$ \Omega_i\le C\frac{t^{kp}}{r_i^{kp}}
\sum_{K'\in\pi'_i} |(2K')\cap S|\, \Ec_k(f;2K')^p_{L_u(S)}.
$$
\par We have assumed that the family of cubes
$\{2K':~K'\in\pi'_i\}$ is a packing (consisting of equal
cubes of diameter $4r_i$). Therefore by definition \rf{OKP}
and by property \rf{QMON} of $\Omega_{k,p}$ we have
$$
\Omega_i\le C\,t^{kp}\,
\frac{\Omega_{k,p}(f;4r_i)^p_{L_u(S)}}{(4r_i)^{kp}}\le
C\,t^{kp}\,\int_{4r_i}^{8r_i} \left
(\frac{\Omega_{k,p}(f;\tau)^p_{L_u(S)}}{\tau^k}\right)^p
\frac{d\tau}{\tau}.
$$
\par Summarizing these estimates for all $i=0,...,m$ we
obtain
$$
I_1:=\sum_{i=0}^m\Omega_i\le
C\,t^{kp}\,\int_{2^{-m+2}\delta_S}^{8\delta_S}
\left
(\frac{\Omega_{k,p}(f;\tau)^p_{L_u(S)}}{\tau^k}\right)^p
\frac{d\tau}{\tau}.
$$
But by \rf{TD} $2^{-m}\delta_S\ge 4000t$ so that
$$
I_1\le C\,t^{kp}\,\int_{t/8}^{8\delta_S} \left
(\frac{\Omega_{k,p}(f;\tau)^p_{L_u(S)}}{\tau^k}\right)^p
\frac{d\tau}{\tau}\,\,, \ \ \ \ \ 0<t\le\delta_S/4000.
$$
\par Let us estimate $\Omega_i$ for $i<0.$ In this case by
\rf{IM} and Theorem \ref{ExtLoc} for every $Q\in\pi_i$,
$i<0$, we have
\bel{NEW} \Ec_k(\tf;Q)_{L_u}
 \le  C\frac{t^k}{\dist(x,S)^k}
\,\Ec_0(f;K^{(x_Q,t)})_{L_u(S)}\,. \ee
\par We continue the proof following the same scheme as for
the case $0\le i\le m$ but using estimate \rf{NEW} rather
than \rf{ST1}. Then we obtain an analog of estimate
\rf{ST2} in the form
$$ \Omega_i:=\sum_{Q\in \pi_i}
|Q|\,\Ec_k(\tf;Q)^p_{L_u}\le C\frac{t^{kp}}{r_i^{kp}}
\sum_{K'\in\pi'_i} |2K'|\, \Ec_0(f;2K')^p_{L_u(S)}. $$
Since $u\le p$, by the H\"{o}lder inequality
$$ |2K'|\, \Ec_0(f;2K')^p_{L_u(S)}=|2K'|
 \left(\frac{1}{|2K'|}\int_{2K'\cap S}|f|^udy\right)^
{\frac{p}{u}} \le \int_{2K'\cap S}|f|^pdy.
$$
Recall that the family $\{2K': K'\in\pi'_i\}$ is a packing
so that
$$
\Omega_i\le C\frac{t^{kp}}{r_i^{kp}} \sum_{K'\in\pi'_i}
\int_{2K'\cap S}|f|^pdy
 \le C\frac{t^{kp}}{r_i^{kp}} \int_S |f|^pdy=
 C\frac{t^{kp}}{r_i^{kp}}\|f\|^p_{L_p(S)}.
$$
Hence
$$
I_2:=\sum_{i<0}\Omega_i\le
C{t^{kp}}\|f\|^p_{L_p(S)}\sum_{i<0}r_i^{-kp}=
C{t^{kp}}\delta_S^{-kp}\|f\|^p_{L_p(S)}\sum_{i<0}2^{ikp}
\le C{t^{kp}}\|f\|^p_{L_p(S)}.
$$
\par Finally, we obtain
$$
\sum_{Q\in \pi} |Q|\,\Ec_k(\tf;Q)^p_{L_u}\le I_1+I_2\le
C\,t^{kp}\,\left\{\int_{t/8}^{8\delta_S} \left
(\frac{\Omega_{k,p}(f;\tau)^p_{L_u(S)}}{\tau^k}\right)^p
\frac{d\tau}{\tau}+\|f\|^p_{L_p(S)}\right\}.
$$
We recall that $\pi$ is an arbitrary packing of equal cubes
with diameter $t/2$ so that by definition \rf{KPU} we
obtain
$$
\Omega_{k,p}(\tf;t/2)_{L_u}\le
C\,t^{k}\,\left\{\left(\int_{t/8}^{8\delta_S} \left
(\frac{\Omega_{k,p}(f;\tau)^p_{L_u(S)}}{\tau^k}\right)^p
\frac{d\tau}{\tau}\right)^{\frac{1}{p}}
+\|f\|_{L_p(S)}\right\}, ~0<t\le \frac{\delta_S}{4000}.
$$
\par To finish the proof of \rf{OM} we observe that for
$t/2\in[\min\{8\delta_S,1/2\},8\delta_S]$ we have
\bel{BT} \Omega_{k,p}(f;t)_{L_u(S)}\le C\|f\|_{L_p(S)}. \ee
(This immediately follows from definition of
$\Omega_{k,p}$\,, see \rf{OKP}, and the H\"{o}lder
inequality.) This allows us to replace the upper limit in
the latter integral by $1/2$ proving that inequality
\rf{OM} is true for $0<t\le \delta_S/8000$.
\par It remains to note that inequality \rf{BT} is true for
$t/2\in [\min\{\delta_S/4000,1/2\},1/2]$ as well which
immediately implies that \rf{OM} is true on all of the
segment $[0,1]$. The proof of inequality \rf{OM} is
finished and we are done.\bx
\begin{remark} \label{MEN}
Inequality \rf{MS} follows from equivalence \rf{EQM} and
Theorem \ref{KPS} with $u=p$.
\end{remark}
\par {P\,r\,o\,o\,f} of Theorem \ref{EXT3}. Clearly,
$\Ec_k(f;Q)_{L_u(S)}\le \Ec_k(F;Q)_{L_u}$
where $F\in L_{u,\,loc}(\RN)$ is an arbitrary extension of
$f$ on all of $\RN$ and $Q$ is an arbitrary cube centered
in $S$. Hence
$$
\|\Ec_k(f;Q(\cdot,t))_{L_u(S)}\|_{L_p(S)}\le
\|\Ec_k(F;Q(\cdot,t))_{L_u}\|_{L_p(\RN)}
$$
so that by \rf{ChB}
$$
I:=\|f\|_{L^p(S)}+ \left(\int_0^1
\left(\frac{\|\Ec_k(f;Q(\cdot,t))_{L_u(S)}\|_{L_p(S)}}
{t^s}\right)^q\, \frac{dt}{t}\right)^{\frac{1}{q}}
 \le C \|f\|_{\BS|_S}.
$$
\par Let us prove the opposite inequality. Using
Theorem \ref{KPS} and the Hardy inequality we obtain
$$
J:=\left(\int_0^1
\left(\frac{\|\Ec_k(\tf;Q(\cdot,t))_{L_u}\|_{L_p(S)}}
{t^s}\right)^q\, \frac{dt}{t}\right)^{\frac{1}{q}}\le C\,I.
$$
We also recall that by Proposition \ref{ELU}
$\|\tf\|_{L_p(\RN)}\le C\|f\|_{L_p(S)}$, so that
$$
\|f\|_{\BS|_S}\le \|\tf\|_{\BS}\approx
\|\tf\|_{L_p(\RN)}+J\le C(\|f\|_{L_p(S)}+I).
$$
\par Theorem \ref{EXT3} is proved.\bx
\begin{remark} The proof of Theorem \ref{KPS}
actually contains the following inequality: for every
$1\le u\le p\le\infty$ and $f\in L_{u,\,loc}(S)$
$$
\Omega_{k,p}(\tf;t)_{L_u}\le
C\,t^{k}\,\left\{\left(\int_{t}^{1} \left
(\frac{\Omega_{k,p}(f;\tau)^p_{L_u(S)}}{\tau^k}\right)^p
\frac{d\tau}{\tau}\right)^{\frac{1}{p}}
+\|f\|_{L_p(S)}\right\}, ~0<t\le 1,
$$
cf. \rf{OM}. This estimate was proved in \cite{S1}, see
also \cite{S2}. Using this inequality rather than the
inequality of Theorem \ref{KPS} one can prove that for
$0<s<k$, $1\le u\le p\le\infty$ and $0< q\le\infty,$
\bel{BSH}
\|f\|_{\BS|_S}\approx \|f\|_{L^p(S)}+
\left(\int_0^1 \left(\frac{\Omega_{k,p}(f;\tau)_{L_u(S)}}
{t^s}\right)^q\, \frac{dt}{t}\right)^{\frac{1}{q}}.
\ee
This version of Theorem \ref{EXT3} has been proved in
\cite{S1}. For the case $1\le p=q\le\infty$ and $s>0$ is
non-integer, description \rf{BSH} was announced in
\cite{Br5}; see also \cite{JW}, p. 211, for another proof
of this result.
\end{remark}
{\bf Acknowledgement.} I am very grateful to Yu. Brudnyi
for stimulating discussions and valuable advice. I also
thank M. Cwikel for helpful comments and remarks.

\end{document}